\documentclass{article}

\usepackage{arxiv}

\usepackage[utf8]{inputenc} 
\usepackage[T1]{fontenc}    
\usepackage{hyperref}       
\usepackage{url}            
\usepackage{booktabs}       
\usepackage{amsfonts}       
\usepackage{nicefrac}       
\usepackage{microtype}      
\usepackage{lipsum}		
\usepackage{graphicx}
\usepackage{doi}
\usepackage{pdfpages}
\usepackage{floatrow}
\usepackage{subcaption}
\usepackage{tikz}
\usepackage{pgfplots}
\usepackage{pgfplotstable}
\pgfplotsset{compat=1.17}
\usepackage{bigints}
\allowdisplaybreaks

\title{Uncountable Infinite Exact Solutions to the FitzHugh-Nagumo Model}

\author{ \href{https://orcid.org/0000-0002-7764-7529}{\includegraphics[scale=0.06]{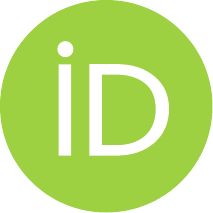}\hspace{1mm}Shahid Sultan Ali Ramji}\thanks{Abdus Salam School of Mathematical Sciences} \\
	Department of Mathematics\\
	University of Karachi\\
	\texttt{shahidsultanali@uok.edu.pk} \\
	\And
	\href{https://orcid.org/0000-0001-5485-4034}{\includegraphics[scale=0.06]{orcid.pdf}\hspace{1mm}Eddy Kwessi} \\
	Trinity University\\
	Texas, U.S.A.\\
	\texttt{} \\
 \And
	\href{https://orcid.org/0000-0001-5528-1207}{\includegraphics[scale=0.06]{orcid.pdf}\hspace{1mm}Mujahid Abbas} \\
	Government College University\\
	Lahore, Pakistan \\
	\texttt{} \\
}

\hypersetup{
pdftitle={Uncountable Infinite Exact Solutions to the FitzHugh-Nagumo Model},
pdfsubject={q-bio.NC, q-bio.QM},
pdfauthor={David S.~Hippocampus, Elias D.~Striatum},
pdfkeywords={First keyword, Second keyword, More},
}

\begin{document}
\maketitle

\begin{abstract}
First time in six decades, uncountable infinite exact solutions of FitzHugh-Nagumo model with diffusion have been found. FitzHugh-Nagumo model is a nonlinear dynamical system applicable to neurosciences, chemical kinetics, cell division, population dynamics, electronics, epidemiology, cardiac physiology and pattern formation. Non-classical symmetry analysis has been carried out to find invariant solutions. In many ways this finding makes the corpus of asymptotics and numerics obsolete. Non singular, physically stable and meaningful  solutions could now be found without distorting the actual model or introducing forced conditions.
\end{abstract}

\section{\label{sec:level1}Introduction:\protect\\}
Embedded within the cells of nervous system,  exists a composition of a special classification termed neurons. This classification of the neurons is a highly significant aspect acting as the main worker units in the nervous system. The neurons are designed in a very unique way in the communication system of the gland cells and muscles thus, enabling them importantly, to act as carriers of information transmission between them.  It is very natural to think that any model for neurons would be complicated. FitzHugh-Nagumo model with diffusion (FHND) arose as a special case of the famed Hodgkin-Huxley (Nobel Prize 1963) model. The pioneering work by Nagumo \cite{Nag62} opened several avenues for research in neuroscience. As a result there exists many different variants of FHND. It can be used to model how voltage of neuron is changing as a function of time and space,
\begin{eqnarray} \label{fd}
\frac{\partial u}{\partial t}&=&D \frac{\partial^{2} u}{\partial x^{2}}-v+g(u) \nonumber\\
\frac{\partial v}{\partial t}&=&\varepsilon[-\beta v+u],
\end{eqnarray}
where $g(u)=u -\frac{u^{3}}{3}$ and $D, \varepsilon,$ and $ \beta$  are real parameters. The holy grail are the integral surfaces $u(x,t)$ and $v(x,t)$. These quantities are known as voltage potential (fast variable) and recovery variable (slow variable) respectively. 
It is a well-established neurodynamic model. For decades this model has validated itself against empirical data. A cursory Google search produces more than 200,000 results for FHN. According to Web of Science FHN has 2564 primary articles and 50930 primary citations after 1970.   Despite of its importance, this original version of FHND has never been cracked exactly before. Although if we either change $g(u)$ slightly or introduce a diffusion term with Laplacian in the second equation or impose some extra conditions on the model then it becomes very easy to work out numerous exact solutions. In the literature, many approaches have been proposed to study its dynamics. In \cite{shah}, we proved that an exact discrete scheme exist for FitzHugh-Nagumo model. We also worked out discrete steady states for our version of discrete FitzHugh-Nagumo model. Later we identified the problem to carry out stability analysis of FitzHugh-Nagumo with diffusion. In the the same article we found numerous steady states for the original continuous FitzHugh-Nagumo model with diffusion. Most of steady states found were space dependent. As pointed out linear stability analysis was not a reliable tool. Luckily we were also able to find isolated fixed points. Here in this manuscript we presenting the Jacobian for FitzHugh-Nagumo model also.  \\
 
Isolated fixed points found in \cite{shah}:  
\begin{eqnarray}\label{fp}
u(t,x)\vspace{0.1cm}&=\frac{\left(\left(4 \sqrt{-\frac{4 \beta^{3} -12 \beta^{2}+12 \beta -4}{\beta}}\right) \beta^{2}\right)^{\frac{1}{3}}}{2 \beta}\nonumber\\ \nonumber\\
&+  \frac{2 (\beta -1)}{\left(\left(4 \sqrt{-\frac{4 \beta^{3}-12 \beta^{2}+12 \beta -4}{\beta}}\right) \beta^{2}\right)^{\frac{1}{3}}}\;, 
\end{eqnarray}
\begin{eqnarray}\label{fp2}
v(t,x)\vspace{0.1cm}&= -\frac{1}{\beta}\left( -\frac{\left(\left(4 \sqrt{-\frac{4 \beta^{3} -12 \beta^{2}+12 \beta -4}{\beta}}\right)\beta^{2}\right)^{\frac{1}{3}}}{2 \beta} \right)\nonumber\\
&+\frac{1}{\beta}\left(\frac{2 (\beta -1)}{\left(\left(4 \sqrt{-\frac{4 \beta^{3} -12 \beta^{2}+12 \beta -4}{\beta}}\right) \beta^{2}\right)^{\frac{1}{3}}} \right)\;.
\end{eqnarray}

\begin{eqnarray}\label{fp3}
u(t,x)=\frac{\sqrt{3}\, \sqrt{\beta-1}}{\sqrt{\beta}}, \ v(t,x)=\frac{\sqrt{3}\, \sqrt{\beta-1}}{\beta^{\frac{3}{2}}}.
\end{eqnarray}
\begin{eqnarray}\label{fp4}
u(t,x)=-\frac{\sqrt{3}\, \sqrt{\beta-1}}{\sqrt{\beta}}, \ v(t,x)=-\frac{\sqrt{3}\, \sqrt{\beta-1}}{\beta^{\frac{3}{2}}}.    
\end{eqnarray}
\begin{eqnarray}\label{fp5}
u(t,x)=0, \ v(t,x)=0.    
\end{eqnarray}
\begin{eqnarray}    \label{fp6}
u(t,x)&=\frac{\left(\sqrt{-\frac{(\beta-1)^{3}}{\beta}}\, \beta^{2}\right)^{\frac{2}{3}}+\beta^{2}-\beta}{\left(\sqrt{-\frac{(\beta-1)^{3}}{\beta}}\, \beta^{2}\right)^{\frac{1}{3}} \beta},\nonumber\\
v(t,x)&=-\frac{-\frac{\left(\sqrt{-\frac{(\beta-1)^{3}}{\beta}}\, \beta^{2}\right)^{\frac{2}{3}}+\beta^{2}-\beta}{\left(\sqrt{-\frac{(\beta-1)^{3}}{\beta}}\, \beta^{2}\right)^{\frac{1}{3}} \beta}}{\beta}.
\end{eqnarray}
Space dependent steady states found in \cite{shah}


\begin{eqnarray}
   u(t,x)=-\frac{\tanh(\frac{\sqrt{2}\, \sqrt{\frac{\beta-1}{\mathrm{D} \beta}}\, (x+\mathit{x0})}{2}) \sqrt{3}\, \sqrt{\beta-1}}{\sqrt{\beta}},\nonumber\\v(t,x)=-\frac{\tanh(\frac{\sqrt{2}\, \sqrt{\frac{\beta-1}{\mathrm{D} \beta}}\, (x+\mathit{x0})}{2}) \sqrt{3}\, \sqrt{\beta-1}}{\beta^{\frac{3}{2}}}. 
\end{eqnarray}
\begin{eqnarray}
u(t,x)=\frac{\tanh(\frac{\sqrt{2}\, \sqrt{\frac{\beta-1}{\mathrm{D} \beta}}\, (x+\mathit{x0})}{2}) \sqrt{3}\, \sqrt{\beta-1}}{\sqrt{\beta}},\nonumber\\v(t,x)=\frac{\tanh(\frac{\sqrt{2}\, \sqrt{\frac{\beta-1}{\mathrm{D} \beta}}\, (x+\mathit{x0})}{2}) \sqrt{3}\, \sqrt{\beta-1}}{\beta^{\frac{3}{2}}}. 
\end{eqnarray}

Using Similariry Solutions a verified group invariant steady state solution was obtained.
\begin{eqnarray}\label{jac fp}
&&u\!\left(x ,t \right)= c_{2} \sqrt{6}\,\sqrt{\frac{\beta -1}{\beta 
  c_{2}^{2}+5 \beta -6}}\
 \nonumber\\ &&\mathrm{JacobiSN}\!\left(\frac{\left(6 c_{1} \mathrm{D} \beta +\sqrt{6}\, \sqrt{\mathrm{D} \beta  \left(-6+5 \beta \right)}\, x \right) \sqrt{6}\, \sqrt{\frac{\beta -1}{\beta  c_{2}^{2}+5 \beta -6}}}{6 \mathrm{D} \beta},\frac{c_{2} \sqrt{5 \beta^{2}-6 \beta}}{-6+5 \beta}\right).\nonumber\\
\end{eqnarray}

There is limited research available on the similarity analysis of FHND, although some articles focus on the conditional symmetry \cite{FMNAG} of the fast equation. Slow variable is always neglected. This research carries out a thorough symmetry analysis for the reaction diffusion system (\ref{fd}). Conditional symmetry analysis has been worked out to find invariant solutions. Yet again this investigation establishes the potential of symmetry methods to unravel non-linearity. First time ever in the sixty two years history of this paradigm infinite exact solutions have been found. Annexure collects some exact results obtained. These results are verified here using the pdetest function available in Maple's library.
\section{Lie Symmetry Analysis}
This part will try to serve as an introduction to lie symmetry analysis. The aim is to use symmetries to solve differential equations. Like any other mathematical persuasion, one needs a definition to proceed. Mostly differential equations are solved by changing the variables involved. 
\begin{equation}
    x^{*}=x^{*}(x,y), \ \ y^{*}= y^{*}(x,y)
\end{equation}
Now there are at least two ways to view this transformation. One can either geometrically view $XY$ plane as the domain of this mapping and $X^{*}Y^{*}$ plane as the codomain. Or one can superimpose $XY$ and $X^{*}Y^{*}$ plane one over another with the common origin. So primarily this mapping maps $(x,y)$ to $(x^{*},y^{*})$. As goal is to built continuous groups these transformation depends on an arbitrary parameter $\epsilon$,
\begin{equation}
    x^{*}= x^{*}(x,y,\epsilon) \ \  y^{*}=y^{*}(x,y,\epsilon).
\end{equation}
Definition of a group requires these elements must locally be  closed under composition, must contain the identity element and invertible to posses inverses. 
\begin{equation}
       x^{**}=x^{**}(x^{*},y^{*},\epsilon^{*})=x^{**}(x,y,\epsilon^{**}) 
\end{equation}
\begin{equation}
    x^{*}(x,y,0)=x; \ \  y^{*}(x,y,0)=y.
\end{equation}
\begin{equation} 
     det \left| \begin{array}{cc}
     \frac{\partial x^{*}}{\partial x} & \frac{\partial x^{*}}{\partial y} \\
      \frac{\partial y^{*}}{\partial x} & \frac{\partial y^{*}}{\partial y} 
      \end{array}\right|\neq 0\;,
   \end{equation}
There are ways to define lie symmetry. A very useful way \cite{hs} is $[X,A]= \lambda A$.  To connect differential equations with the operator $A$ consider an example. Let $y=(x,a^\alpha),\alpha= 1,\dots,n$ be the general solution of a differential equation. For the corresponding partial differential equation $Af=0$, independent solutions $\phi^\alpha, \alpha= 1,\dots,n$ can be obtained by differentiation of the solution set. Differentiation of solution system $n$ times results in a new system of equations. Generalized version of implicit function theorem guarantees that $a^{\alpha}$ can be made the subject of this system. That is to say $a^{\alpha}$ can be written explicitly and appearing left side of the system. These $n$ constants are precisely the level surfaces of $n$ linearly independent solutions to the associated pde $Af=0,$ related to the given $nth$ order ode.  In particular the general solution of $y^{\prime \prime}+y=0$ is given by
\begin{equation}\label{ls1}
y=\left(\phi_0^1\right)^{1 / 2} \sin \left(x-\phi_0^2\right). 
\end{equation}
\begin{equation}\label{ls2}
\text{So} \  y^{\prime}=\left(\phi_0^1\right)^{1 / 2} \cos \left(x-\phi_0^2\right),    
\end{equation}
\begin{equation}
y^{\prime \prime}=-\left(\phi_0^1\right)^{1 / 2} \sin \left(x-\phi_0^2\right).       
\end{equation}
Dividing (\ref{ls1}) by (\ref{ls2})
$$\tan \left(x-\phi_0^2\right)=y / y^{\prime}, $$
$$
x-\phi_0^2 = \tan ^{-1}\left(y / y^{\prime}\right),$$
\begin{equation}\label{phi2}
\phi_0^2=x-\tan ^{-1} y / y^{\prime}.    \end{equation}
Squaring and adding (\ref{ls1}) and (\ref{ls2})
\begin{equation}\label{phi1}
y^2+y^{\prime^2}=\phi_0^1.    
\end{equation}
(\ref{phi2}) and (\ref{phi1}) are the two linearly independent solutions of $$(\frac{\partial}{\partial x}+y^{'} \frac{\partial}{\partial y}-y \frac{\partial}{\partial y^{'}})f=0.$$ In fact they are functionally independent. \ \ $\triangleleft$ \\

\ \ Given $X$, a lie symmetry of $A$ i.e. $[X,A]=\lambda A$. Define $$\hat{X}=X +\mu(x,y,y^{\prime},\dots,y^{n-1})A$$ then for some $\hat{\lambda},$ $[\hat{X}, A]=\hat{\lambda}A$.   
\begin{eqnarray*}
[\hat{X}, A]&=&[X+\mu A, A] \\
& =&[X, A]+[\mu A, A] \\
& =&\lambda A+(\mu A) A - A(\mu A) \\
& =&\lambda A+(\mu A) A-(A \mu) A-\mu(A(A) \\
& =&(\lambda-A \mu) A \\
& =&\hat{\lambda} A 
\end{eqnarray*}

The question is when $\hat{X}$ is going to be a lie symmetry for $A$? As already worked in \cite{hs}, $\lambda=-A \xi$ and  $\xi$ is a function $x$ and $y$. In order for $\hat{X}$ to be a symmetry $\mu$ should be independent of the derivatives $y^{\text {(i)}}$. On the contrary let us say $\mu$ has the term $y^{(n-1)}$ then the presence of
$\frac{\partial}{\partial x}$ in the total differential $A$ will throw the $\hat{\lambda}$ into the world of $y^{(n)}$ and hence $\hat{X}$ will not remain a symmetry. \ \ $\triangleleft$\\

In order to use lie symmetry understanding of this operator $A$ could be very useful. Let us explore more about it. $\hat{A}=\mu(x,y,\dots,y^{n-1})A$ and A are equivalent as their solutions $\phi^{a}$ are same. Let us show that for every symmetry $X$ of $A$, a commutator $\hat{A}$ with $X$ exists. $\mu$ is already a yanked out factor of $A$ so obviously anything that satisfies $A$ necessarily satisfies $\mu A$. Proof of the statement resides in the following argument. 
\begin{eqnarray*}
[X, \hat{A}]&=&[X, \mu A]=X(\mu A)-(\mu A) X \\
&=&(X \mu) A+\mu(X A)-\mu(A X) \\
&=&(X \mu) A+\mu(X A-A X) \\
& =&(X \mu) A+\mu \lambda A \\
& =&(X \mu+\mu \lambda) A
\end{eqnarray*}

Now $X \mu+\mu \lambda=0$ is a pde whose solution space is non empty, hence for every lie symmetry commutator exists. \ \ $\triangleleft$\\

  \ \ As pointed out earlier pragmatism will be the philosophy for this discourse. Hence the best approach would be exemplification. So let us directly integrate the following first order ordinary differential equations by their respective lie symmetries. 
\begin{eqnarray*}
 y^{\prime}&=&\frac{-(y+2 x)}{x}, \ \ X=\frac{1}{y+2 x} \frac{\partial}{\partial x} \\
 y^{\prime}&=& \ \ f(y) g(x), \ \ X=\frac{1}{g(x)} \frac{\partial}{\partial x}
\end{eqnarray*}

For the former the first integral can directly be computed.
\begin{eqnarray*}
\phi&=&\int \frac{d y+\frac{(y+2 x)}{x} d x}{-\left(\frac{1}{y+2 x}\right)(\frac{-(y+2x)}{x})} \\
&=&\int x d x+\left(\frac{(y+2 x)}{x}\right)x dx \\
& =&x y+x^2 \\
 y&=&\frac{c-x^2}{x}     
\end{eqnarray*}

As for the second equation, exploiting given symmetry as an integrating factor, conserved quantity, the first integral can directly be given via the following integral.  
\begin{eqnarray*}
\phi&=&\int \frac{d y-f(y)(g(x))dx}{\frac{1}{g(x)} f(y) g(x)} \\
\phi&=&-\int \frac{d y}{f(y)}+\int g(x) d x 
\end{eqnarray*}$\quad\triangleleft$\\ 

\ \ Beside finding solutions of different equations, symmetries could be of great help for order reduction. As an example see 

 $$2 y^{\prime} y^{\prime \prime \prime}-3 y^{\prime \prime^2}=0$$
 admits the following symmetry
$$
X=a \frac{\partial}{\partial x}+b \frac{\partial}{\partial y}+cy \frac{\partial}{\partial y}.
$$
It is not difficult to see that differential equation  does not depend on $y$ so the given equation has the structure
$$
s^n=\hat{\omega}\left(t, s^{\prime}, s^{\prime \prime}, \dots, s^{n-1}\right) .
$$
Hence it is already an equation of order $n-1$ for $y^{\prime}$.\ \ $\triangleleft$\\

\ \ For the equation
$y^{\prime \prime \prime}-y y^{\prime \prime}+y^{\prime 2}=0;$
$$
X=a \frac{\partial}{\partial x}+b\left(x \frac{\partial}{\partial x}-y \frac{\partial}{\partial y}\right)
$$
is a lie point symmetry. If we use translation invariance along $x$ part of this two parameter local group. Then we can reduce one order in the following way.
\begin{eqnarray*}
 X=\frac{\partial}{\partial x}. \\
 t=y, \ \ s=x. \\
 y^{\prime}=\left(\frac{ds}{dt}\right)^{-1}, \\
 y^{\prime \prime}=\frac{-s^{\prime \prime}}{\left(s^{\prime}\right)^3}, \\
 \frac{d^3 s}{d t^3}=\frac{1}{t_{n}+t_{y} y^{\prime}} \frac{d\left(\frac{d s^2}{d r^{2}}\right)}{d x} \\
 =\frac{1}{y^{\prime}} \frac{d\left(\frac{-y^{\prime \prime}}{y^{\prime^3}}\right)}{d x}=\frac{-y^{\prime^3} y^{\prime \prime \prime}+3 y^{\prime\prime^2} y^{\prime^2}}{y^{\prime^7}} \\
 s^{\prime \prime \prime}=\frac{-y^{\prime} y^{\prime \prime \prime}+3 y^{\prime \prime^2}}{y^{\prime^5}}, \\
 y^{\prime \prime \prime}=\frac{s^{\prime \prime \prime}}{s^{\prime^4}}+\frac{3 s^{\prime \prime}}{s^{\prime ^2}}, \\
 \frac{s^{\prime \prime \prime}}{s^{\prime^5}}-\frac{3 s^{\prime \prime^{2}}}{s^{\prime^{6}}}=-y^{\prime } y^{\prime \prime \prime}, \\
 y^{\prime \prime}=\frac{3 s^{\prime \prime^ 2}}{s^{\prime ^5}}-\frac{s^{\prime \prime \prime}}{s^{\prime^4}} = \frac{3 s^{{\prime \prime^2}}-s^{\prime}s^{\prime \prime}}{s^{{\prime^5}}}, \\
\frac{3 s^{\prime \prime^ 2}-s^{ \prime} s^{\prime \prime}}{s^{\prime^5}}- t\left(\frac{-s^{\prime \prime}}{s^{\prime^3}}\right)+\frac{1}{s^{\prime}}=0, \\
3 s^{\prime \prime^2}-s^{\prime} s^{\prime \prime \prime}+s^{\prime \prime} s^{\prime^2} t+s^{\prime^3}=0. 
\end{eqnarray*}\ \ $\triangleleft$

\ \ As another example not for order reduction but to actually find a solution. Transform the following symmetry
$$
X=[a+b(x-y)] \frac{\partial}{\partial x}+a \frac{\partial}{\partial y}
$$
with $b=0$ to coordinates $t, s$ given in
$$
t=y, \ \ s=\ln (x-y), 
$$
and then use the generator to integrate
$$ s^{\prime\prime}+s^{\prime 2}+1=0.$$
\begin{eqnarray*}
 X&=&a \frac{\partial}{\partial x}+a \frac{\partial}{\partial y}. \\
 X&= &X(t) \frac{\partial}{\partial t}+X(s) \frac{\partial}{\partial s}\\
& =&a \frac{\partial}{\partial t}
\end{eqnarray*}

\begin{eqnarray*}
\phi&=&\int \frac{d s^{\prime}-\left(-1-s^{\prime^2}\right) d t}{-a\left(-1-s^{\prime^2}\right)} \\
&=&\frac{1}{a} \int \frac{d s^{\prime}}{1+s^2}+\frac{1}{a} \int d t \\
&=&\frac{1}{a} \tan ^{-1} s^{\prime}+\frac{1}{a} t, \\
 s^{\prime}&=&\tan (a\ \phi-t), \\
 s&=&\int \tan (a\ \phi-t)+C 
\end{eqnarray*}$\quad\triangleleft$





\ \ When the prolonged generators are intransitive in the solution space of differential equation, the integration process fails. Even with these generators, the action of whose group doesn't cover all solution space, differential equations could be solved merely by algebraic manipulation. When the action of these generators divide the solution space into families of different orbits a linear relationship $\mu_1 X_1+\mu_2 X_2$ $+ \nu A=0$ between these operators exists. When $\mu_1 X_1+\mu_2 X_2$ $+ \nu A=0$ holds then $\nu$, $\mu_1$ and $\mu_2$  can't be constants. 
Let us assume that such a relationship  holds with constants, then one of it being nonzero makes, say
\begin{equation} \label{lin rel}
 X_{2}= \Psi X_{1}+ \nu A. \\    
\end{equation}
Here $\Psi$ and $\nu$  are constants. 
\begin{eqnarray*}
 &\left[A, X_{2}\right]&=-\lambda_2 A=( A \Psi) X_{1}-\Psi \lambda_1 A+(A v)A \\
& \Rightarrow& \ \ \lambda_2=\Psi \lambda_1 \\
& \Rightarrow&\left[A, X_2\right]=\Psi\left[A, X_1\right] \\
& \Rightarrow& \ \ X_2=\Psi X_1. \\
\end{eqnarray*}

Putting $X_{2}$ in (\ref{lin rel}), 
$$
 \nu A =0.
$$
But this is nonsensical. Similarly taking only one as constant it can be shown that if such a linear relation hold then $\nu$, $\mu_1$ and $\mu_2$  can't be constants.\ \ $\triangleleft$\\

\ \  This discussion about connecting operators become unavoidable when one wants to solve a differential equation with $G_3IX$. For the time being revert to integration process again. If  $X_1, X_2$  form a lie algebra, $\delta_1=\xi_1 n_1-\xi_2 n_1=0$ holds then show that $y^{\prime \prime}=\omega\left(x, y, y^{\prime}\right)$ is transformable into $d^2 \hat{s} / d t^2=0$,  admitting eight symmetries.\\

Case I.
\begin{eqnarray*}
{\left[X_1, X_2\right]=0, \ \ \delta=0}. \\
s^{\prime \prime}=\hat{\omega}(t).
\end{eqnarray*}

Lets try
\begin{eqnarray*}
 \hat{s}(t)&=&s-\int\int^t \hat{\omega}\left(t^{\prime}\right) d t^{\prime} d t+\phi_0 t+\psi_0, \\
 \frac{d \hat{s}(t)}{d t}&=&\frac{d s}{d t}-\int \hat{\omega}(t) d t+\phi_0, \\
 \frac{\left.d^2 \hat{s} (t\right)}{d t^2}&=&\frac{d^2 s}{d t^2}-\hat{\omega}(t)=0.
\end{eqnarray*}

From Chapter 4 Question 2 \cite{hs}, the last equation is transformable to trivial second order differential equation and they both admit eight symmetries.\\

Case II.
$$
\begin{array}{r}
{\left[X_1, X_2\right]=X_1, \ \ \delta=0}. \\

s^{\prime\prime}=s^{\prime} \hat{\omega}(t).
\end{array}
$$
Lets try 
\begin{eqnarray*}
\hat{s}(t)=s-\int \phi_{0} e^{\hat{\omega} (t^{\prime}) d t^{\prime}} d t+\phi_{0}, \\
\frac{d \hat{s}(t)}{d t}=\frac{d s}{d t}-\phi_{0} e^{\int \hat{\omega}(t) d t}, \\
 \frac{d^2 \hat{s}(t)}{d t^2}=\frac{d^2 s}{d t^2}-\phi_{0} e^{\hat{\omega} (t) dt} \cdot \hat{\omega}(t), \\
 \frac{d^{2} s}{d t^{2}}-s^{\prime} \hat{w}(t)=0.
\end{eqnarray*}

Again this second order homogeneous linear differential equation admits eight symmetries as it can be transformed to $y^{\prime\prime}=0$.\ \ $\triangleleft$\\

\ \ When a second order ordinary differential equation admits $G_2$, two independent first integrals can directly be given via line integrals. Hence the solution can be found by eliminiting $y^\prime$ from these two first integrals. Below is an example with an abelian lie algebra.\\

Solve $y^{\prime \prime}=y+x^2$ using 
$X_{1}=e^x\frac{\partial}{\partial y}, X_{2}=e^{-x}\frac{\partial}{\partial y}$. \\
$${\left[X_1, X_2\right] }  =\left(e^x\left(\frac{\partial}{\partial y}\left(e^{-x}\right)-e^{-x} \frac{\partial}{\partial y}\left(e^x\right)\right) \frac{\partial}{\partial y}\right) =0.$$

\begin{eqnarray*}  
\delta&=&0. \ \\
 \eta_{1}^{\prime}&=&\frac{\partial \eta_{1}}{\partial x}+\eta^{\prime}\left(\frac{\partial \eta_{1}}{\partial_y}-\frac{\partial \xi_{1}}{\partial x}\right)-{y}^{\prime^{2} } \frac{\partial \xi_{1}}{\partial y} =e^x. \\ 
\eta_{2}^{\prime}&=&-e^{-x}. \\
\end{eqnarray*}
\begin{eqnarray*}
\phi&=& \bigintsss\frac{\left| \begin{array}{ccc}
d x & d y & d y^{\prime} \\
1 & y^{\prime} & y+x^{2} \\
0 & e^x & e^x
\end{array}\right|}{\left|\begin{array}{ccc}
1 & y^{\prime} &  y+x \\
0 & e^{x} & e^{x} \\
0 & e^{-x} & -e^{-x}
\end{array}\right|} \\
& =&\frac{-1}{2} \int\left(y^{\prime} e^x-y e^x-x^2 e^x\right) d x+e^x d y+e^x d y^{\prime} \\
& =&-\frac{1}{2} \int d\left(y^{\prime} e^x\right)-d\left(y e^x\right)- x^2 e^x d x \\
& =&-\frac{1}{2}\left(e^x y^{\prime}\right)+\frac{1}{2}\left(e^x y\right)+\frac{1}{2} \int x^2 e^x d x \\
& =&\frac{1}{2} e^x\left(y-y^{\prime}\right)+\frac{1}{2}\left[x^2 e^x-\int 2 x e^x d x\right] \\
& =&\frac{1}{2} e^x\left(y-y^{\prime}\right)+\frac{1}{2}\left[x^2 e^x-2 x e^x+2 e^x\right] \\
& \phi=&\frac{1}{2} e^x\left(y-y^{\prime}+x^2+2 x+2\right) \\
& \psi=&-\frac{1}{2} \int\left|\begin{array}{ccc}
d x & d y &  { dy^{\prime} } \\
1 & y^{\prime} & y+x^2 \\
0 & e^{-x} & -e^{-x}
\end{array}\right| \\
& =&\frac{-1}{2} \int\left(-e^{-x} y^{\prime}-y e^{-x}-x^2 e^{-x}\right) d x+e^{-x} d y  +e^{-x} dy^{\prime} \\
& =&-\frac{1}{2} \int d\left(y^{\prime} e^{-x}\right)+d\left(y e^{-x}\right)-x^2 e^{-x} d x \\
& =&-\frac{1}{2} e^{-x} y^{\prime}-\frac{1}{2} y e^{-x}+\frac{1}{2} \int x^{2} e^{-x} d x \\
& =&-\frac{1}{2} e^{-x}\left(y+y^{\prime}\right)+\frac{1}{2}\left(-x^2 e^{-x}-2 x e^{-x}-2 e^{-x}\right) \\
 \psi&=&\frac{-1}{2} e^{-x}\left(y+y^{\prime}+x^{2}+2 x+2\right), \\
 2 e^{-x} \psi &=&  y-y^{\prime}+x^2-2 x+2, \\
-2 e^x \phi &=&y+y^{\prime}+x^2+2 x+2, \\
 y&=&\phi e^{-x}-\psi e^x-x^2-2. \\
\end{eqnarray*}

Although Lie's integration procedure has led us to the correct solution but observe that the given differential equation is second order non-homogeneous linear. And simpler and faster methods are already available to find the solution. The problem with the theory of differential equations is that it makes mathematics look more of an art than science. Diverse methods each requiring  ingenuity of the solver if not guesses. Lie's method could unifies all techniques. It has the potential to be 'one size fits all' or the 'theory of everything' for differential equations.\ \ $\triangleleft$ \\ 






\ \ To see lie's method of integration using normal form work, consider  $3 y y ^{\prime\prime}=5y^{\prime^2}$. It admits $X_1=\frac{\partial}{\partial x}$, $X_2=x \frac{\partial}{\partial x}, X_3=y\frac{\partial}{\partial y}$. Its algorithmic to integrate given differential equation using three symmetries. The given equation has no direct dependence on $x$ hence $X_1$, is a symmetry. In fact its not hard to see $X_2$ is also a symmetry. For the sake of completeness, 
\begin{eqnarray*}
 \tilde{x}&=&\epsilon x, \\
 \tilde{y}&=&y. \\
 \frac{d \tilde{y}}{d \tilde{x}}=\frac{d y}{d(\epsilon x)} &=& \frac{1}{\epsilon} \frac{d y}{d x}, \\
 \frac{d^2 \tilde{y}}{d \tilde{x}^{2}}&=&\frac{d}{d \tilde{x}}\left(\frac{1}{\epsilon} \frac{d y}{d x}\right)=\frac{1}{\epsilon^{2}} \frac{d ^{2}y}{d x^2}, \\
3 y \frac{1}{\epsilon^2} \frac{d ^{2} y}{d x^2}&=&5\left(\frac{1}{\epsilon} \frac{d y}{d x}\right)^2 \\
 \Rightarrow 3 y y^{\prime\prime}&=&5 y^{\prime^2}; \\
\end{eqnarray*}

$X_2$ is a symmetry. Now take
\begin{eqnarray*}
 \tilde{x}&=&x, \\
\tilde{y}&=&\epsilon y. \\
 \frac{d \tilde{y}}{d \tilde{x}}&=&\epsilon\frac{d y}{d x}, \\
 \frac{d^2 \tilde{y}}{d \tilde{x}^2}&=&\epsilon\frac{d^2 y}{d x^2}, \\
 3 \epsilon y \epsilon y^{\prime \prime}&=&5 \epsilon^{2} y^{\prime 2}, \\
\Rightarrow 3 y y^{\prime \prime}&=&5 y^{\prime^{2}}; \\
\end{eqnarray*}

$X_3$ is a symmetry.
Now take $X_1=\frac{\partial}{\partial x},  X_3=y \frac{\partial}{\partial y}$.
These are commuting operators and for them, 
$$
\Delta=\left| \begin{array}{ccc}
1 & y^{\prime} & \frac{5 y^{\prime^2}}{3 y} \\
1 & 0 & 0 \\
0 & y & y^{\prime} \\
\end{array}\right|= \frac{2 y^{\prime^2}}{3}.
$$

Nonzero determinant means the lie algebra generated by $X_1$ and $X_3$ acts transitively in the space of first integrals. $ \phi, \psi$ are now guaranteed.

\begin{eqnarray*}
  \phi&=&\bigintsss  \frac{\left|\begin{array}{ccc}
d x & d y & d y^{\prime} \\
1 & y^{\prime} & \frac{5 y^{\prime^2}}{3 y} \\
1 & 0 & 0
\end{array} \right|\,}{\frac{2 y^{\prime^2}}{3}} \\
& =&\int \frac{\frac{5 y^{\prime^2}}{3 y} d y-y^{\prime} d y}{\frac{2 y^{\prime^2}}{3}}, \\
 \phi&=&\frac{5}{2} \ln |y|-\frac{3}{2} \ln \left|y^{\prime}\right|. \\
 \phi_{0}&=&\frac{y^{\prime}}{y^{5 / 3}}. \\
 \psi&=&\bigintsss \frac{\left|\begin{array}{ccc}
d x & d y & d y^{\prime} \\
1 & y^{\prime} & \frac{5 y^{\prime}}{3 y} \\
0 & y & y^{\prime}
\end{array} \right|\,}{\frac{2 y^{\prime^2}}{3}} \\
& =&\int \frac{3}{2 y^{\prime^2}}\left[\left(y^{\prime^2}-\frac{5 y^{\prime^2}}{3}\right) d x-y^{\prime} d y+y d y^{\prime}\right] \\
& =&\frac{3}{2} y^{\prime^2}\left(-\frac{2}{3} 5 y^{\prime^2} x\right)-\frac{3}{2 y^{\prime ^2}} y^{\prime} y-\frac{3 y}{2 y^{\prime}}, \\
 \psi_{0}&=& -5 x-\frac{3 y}{y^{\prime}}, \\
-5 x y^{\prime}-3 y&=&y^{\prime} \psi_{0}^{}, \\
 y^{\prime}&=&\frac{3 y}{\psi_{0} +5 x}, \\
\phi_{0} y^{5 / 3}&=&\frac{3 y}{\psi_{0}+5 x}, \\
y^{2 / 3}&=&\frac{3}{\phi_{0} \left(\psi_{0}+5 x\right)}, \\
 y&=&\frac{c_1}{\left(c_{2}+5 x\right)^{3 / 2}}.
\end{eqnarray*}\ \ $\triangleleft$

\ \ In order to solve differential equations almost always some sort of transformations are involved. At the background of Lie's method are finite transformations, these are local group with parameters. All the information about these finite transformations are contained in infinitisimal generators. Next observe that $G_3 IX$ 

\begin{eqnarray*}
&& \left[X_1, X_2\right]=X_3,  \left[X_2, X_3\right]=X_1, \\
&& \left[X_2, X_1\right]=X_2, 
\end{eqnarray*}

lives inside the projective group 
$$
\begin{aligned}
X & =\left(a_1+a_2 x+a_3 y+a_4 x y+a_5 x^2\right) \frac{\partial}{\partial x} \\
& +\left(a_6+a_7 x+a_8 y+a_5 x y+a_4 y^2\right) \frac{\partial}{\partial y}
\end{aligned}
$$
as a subgroup. Before proving this assertion recall that $G_3 IX$ can be realized in $2-\operatorname{space}$ as
$$
\begin{aligned}
& X_1=\frac{\partial}{\partial x}, \ \ X_2=\cot y \cos x \frac{\partial}{\partial x}, \\
& X_3=-\sin x \cot y  \frac{\partial}{\partial x}+\cos x \frac{\partial}{\partial y}.
\end{aligned}
$$
The task could be accomplished if inside the projective group a subgroup that walks and talks like this representation could be found. The structural constants here when viewed as an abstract lie algebra is known as $so(3)$,(Gothic Font). $so(3)$ generates $SO(3)$. $so(3)$ can be realized as
$$
\begin{aligned}
& X_1=y \partial_{x} -x \partial_{y}, \\
& X_2=\frac{1}{2}\left(1+x^2-y^2\right) \partial_x+x y \partial_{y}, \\
& X_3=x y \partial_{x}+\frac{1}{2}\left(1-x^2+y^2\right) \partial y
\end{aligned}
$$
which loosely speaking shows its presence inside projective group. One can write the commutator table for projective group and look for a block of this three generators. This is a laborious job writing again and again in different basis. Standing on Stepanhi's shoulders \cite{hs} we can write 
$$
\begin{aligned}
& X_1=X_1 \phi \frac{\partial}{\partial \phi}+ X_1 \psi \frac{\partial}{\partial \psi}, \\
& X_2=X_2 \phi \frac{\partial}{\partial \phi} + X_2 \psi \frac{\partial}{\partial \psi}, \\
& X_3=X_3 \phi \frac{\partial}{\partial \phi} + X_3 \psi \frac{\partial}{\partial \psi},
\end{aligned}
$$
as
$$
\begin{aligned}
& X_1= y \frac{\partial}{\partial x}+\left(-x\right) \frac{\partial}{\partial y}, \\
& X_2= -xy \frac{\partial}{\partial x}+\left(-1-y^2\right) \frac{\partial}{\partial y}, \\
& X_3=(1+x^2) \frac{\partial}{\partial x} +x y \frac{\partial}{\partial y}. \\
&
\end{aligned}
$$
Coming from solution space into the $xy$ plane won't change the commutators. The infinitesimal generators of the lie algebra corresponding is projective group

\begin{eqnarray*}
&& \hat{X}_1=x^2 \frac{\partial}{\partial x}+x y \frac{\partial}{\partial y}, \ \ \hat{X}_2=xy \frac{\partial}{\partial x}+y^2 \frac{\partial}{\partial y}, \\
&& \hat{X}_3=x \frac{\partial}{\partial x}, \ \ \hat{X}_4=y \frac{\partial}{\partial x}, \ \ \hat{X}_5=\frac{\partial}{\partial x},\\
&& \hat{X}_6= x\frac{\partial}{\partial y}, \ \ \hat{X}_7=y \frac{\partial}{\partial y}, \ \ \hat{X}_8= \frac{\partial}{\partial y}. \\ \\
&& X_1=\hat{X}_4-\hat{X}_6, \\
&& X_2=-\hat{X}_8-\hat{X_2}, \\
&& X_3=\hat{X_1}+\hat{X_5}.
\end{eqnarray*}

Hence $G3IX$ is a subgroup of the projective group.\ \ $\triangleleft$ \\

\ \ Coming back again to the basic goal of solving differential equations. Two examples earlier showed that presence of $G_2$ is enough to solve a second order ordinary differential equation using line integral. A second order ordinary differential equation can admit up to eight symmetries. So for classification the only problem is any algebra that does not have $G_2$ inside. Only such algebra in the plane is group of rotations. Actually classification of differential equations was a project dear to Lie's heart. Although not for its own sake but for solutions of differential equations.\\

$y^{\prime \prime}=2\left(x y^{\prime}-y\right) / x^2$ admits 
$$
X_1=y \frac{\partial}{\partial x}-x \frac{\partial}{\partial y}, \ \
X_2=\frac{-y}{x}  \frac{\partial}{\partial x}+\left(x^2-\frac{y^2}{x^2}\right) \frac{\partial}{\partial y},
$$
$$
X_3=\left(1+x^2\right) \frac{\partial}{\partial x}+\left(2 x y+\frac{y}{x}\right) \frac{\partial}{\partial y}.$$
It would be a different problem if inside this group there is no two dimensional subgroup with generators  $Y_\alpha=a_\alpha^a X_i, \alpha=1,2$. To check this easiest way is to workout the commutator table and see if this is rotation group or not. Lets take a longer route. 
If one solution $U(x)$ of $\sum_{i=0}^n f_i(x) y^{(i)}=0$ is known then $X=  U(x) \frac{\partial}{\partial x}$ is a symmetry.  $U(x)=x$ is already known

\begin{eqnarray*}
 \text {check:}\ &&y=x, \\
&& y^{\prime}=1, \\
&& y^{\prime \prime}=0, \\
&& x^{2} y^{\prime \prime}-2 x y^{\prime}+2 y=0, \\
&& x^2(0)-2 x(1)+2(x)=0, \\
&& 0=0. \\
\end{eqnarray*}

Symmetry has suggested to introduce $t, s$, so that $X s=1, \ X t=0$.
\begin{eqnarray*}
x \frac{\partial t}{\partial y}=0, & x \frac{\partial s}{\partial y} =1 \\
\Rightarrow t=x, & s=\frac{y}{x}.
\end{eqnarray*}

Prolonging these coordinates
\begin{eqnarray*}
\frac{d s}{d t}&= & \frac{s_x+s_y y^{\prime}}{t_x+t_y y^{\prime}}, \\
s^{\prime}&= & -\frac{y}{x^2}+\frac{y^{\prime}}{x}, \\
s^{\prime \prime}&= & \frac{\left(\frac{\partial}{\partial x}+y^{\prime} \frac{\partial}{\partial y}+y^{\prime \prime} \frac{\partial}{\partial y^{\prime}}\right)\left(-\frac{y}{x^2}+\frac{y^{\prime}}{x}\right)}{\left(\frac{\partial}{\partial x}+y^{\prime} \frac{\partial}{\partial y}+y^{\prime \prime} \frac{\partial}{\partial y^{\prime}}\right) x}\\
& =&\frac{2 y}{x^3}-\frac{y^{\prime}}{x^2}-\frac{y^{\prime}}{x^2}+\frac{y^{\prime \prime}}{x} \\
& =&\frac{2 y}{x^3}-\frac{2 y^{\prime}}{x^2}+\frac{2\left(x y^{\prime}-y\right)}{x^3} \\
& =&\frac{2 y}{x^3}-\frac{2 y^{\prime}}{x^2}+\frac{2 y^{\prime}}{x^2}-\frac{2 y}{x^3}\\
&=&0 ,\\
 s^{\prime \prime}&=&0, \\
 s&=&a t+b, \\
 \frac{y}{x}&=&a x+b, \\
 y&=&a x^2+b x. \\
\end{eqnarray*}

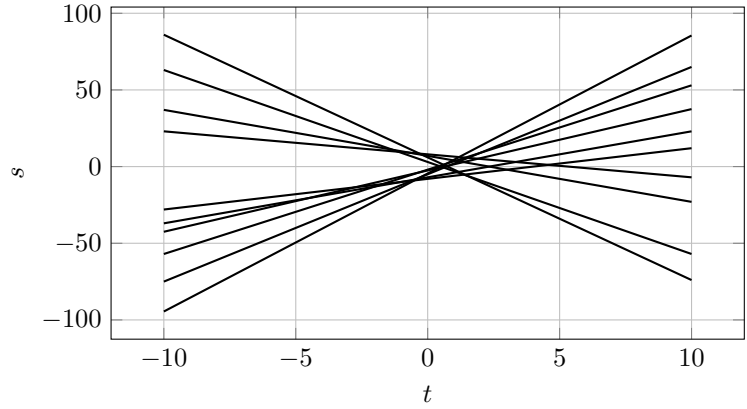
\begin{figure}[ht]
    \centering
    \begin{tikzpicture}
        \begin{axis}[
            title={},
            xlabel={$t$},
            ylabel={$s$},
            grid=major,
            width=10cm,
            height=6cm,
            domain=-10:10
        ]
        \addplot[thick, samples=100] {7*x - 5};
        \addplot[thick, samples=100] {-6*x + 3};
        \addplot[thick, samples=100] {4*x - 2.5};
        \addplot[thick, samples=100] {-1.5*x + 8};
        \addplot[thick, samples=100] {3*x - 7};
        \addplot[thick, samples=100] {9*x - 4.5};
        \addplot[thick, samples=100] {-8*x + 6};
        \addplot[thick, samples=100] {5.5*x - 2};
        \addplot[thick, samples=100] {-3*x + 7};
        \addplot[thick, samples=100] {2*x - 8};
        \end{axis}
    \end{tikzpicture}
    \caption{$st$ plane with random straight line solutions}
    \label{fig:st-plane}
\end{figure}

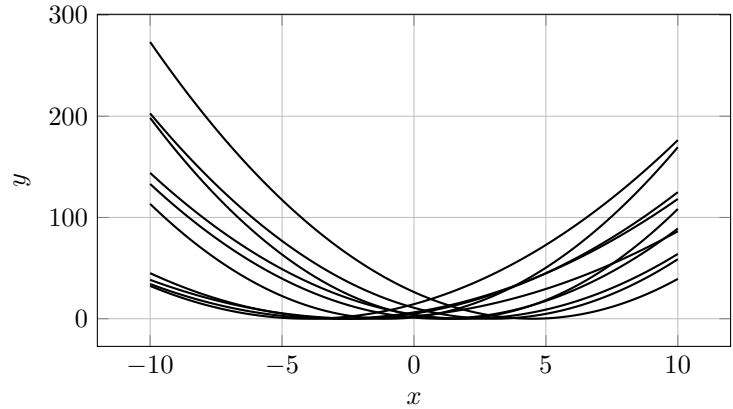
\begin{figure}[ht]
    \centering
    \begin{tikzpicture}
        \begin{axis}[
            title={},
            xlabel={$x$},
            ylabel={$y$},
            grid=major,
            width=10cm,
            height=6cm,
            domain=-10:10
        ]
        \addplot[thick, samples=100] {1.2*(x-3)^2};
        \addplot[thick, samples=100] {0.8*(x+2.5)^2};
        \addplot[thick, samples=100] {1.5*(x-1.5)^2};
        \addplot[thick, samples=100] {0.9*(x+4)^2};
        \addplot[thick, samples=100] {1.3*(x-4.5)^2};
        \addplot[thick, samples=100] {0.7*(x+3)^2};
        \addplot[thick, samples=100] {1*(x-2)^2};
        \addplot[thick, samples=100] {1.4*(x+1)^2};
        \addplot[thick, samples=100] {1.1*(x-1)^2};
        \addplot[thick, samples=100] {0.6*(x+2)^2};
        \end{axis}
    \end{tikzpicture}
    \caption{$xy$ plane with translated parabolas}
    \label{fig:xy-plane}
\end{figure}

The transformation $(x, y) \longrightarrow(t, s)$ has straightened the frame of differential equation. Symmetry unraveled that given differential equation is trivial $ s^{\prime \prime}=0$. The relation  $$X_1=\phi(x,y,y^\prime) X_2+\psi(x,y,y^\prime) X_3+v(x,y,y^\prime) A$$ could have directly given two first integrals if given symmetries were forming $G_3IX$. 
\begin{eqnarray*}
&& {\left[X_1, X_2\right]=\left\{\left(y \frac{\partial}{\partial x}-x \frac{\partial}{\partial y}\right)\left(-\frac{y}{x}\right)-\left(\left(\frac{-y}{x}\right) \frac{\partial}{\partial x}+\left(x^2-\frac{y^2}{x^2}\right) \frac{\partial y}{\partial y}\right)\right\}} \\
&& \frac{\partial}{\partial x}+\{\left(y \frac{\partial}{\partial x}-x \frac{\partial}{\partial y}\right)\left(x^2-\frac{y^2}{x^2}\right)-\left(-\frac{y}{x}\right) \frac{\partial}{\partial x}+\left(x^2-\frac{y^2}{x^2}\right)  \\
&& \left.\frac{\partial}{\partial y}(-x)\right\} \frac{\partial}{\partial y} \\
&& =\left(1-x^2+\frac{2 y^2}{x^2}\right) \frac{\partial}{\partial x}+\left(2 x y+2 \frac{y^3}{x^3}+\frac{3 y}{x}\right) \frac{\partial}{\partial y} \\
&& \neq X_3. \\
\end{eqnarray*}

Changing the coordinates
\begin{eqnarray*}
X_1 (t) & =&\left(y \frac{\partial}{\partial x}-x \frac{\partial}{\partial y}\right) x=y=s t, \\
X_1(t) & =&s t, \\
X_1(s) & =&\left(y \frac{\partial}{\partial x}-x \frac{\partial}{\partial x}\right)\left(\frac{y}{x}\right)=-\frac{y^2}{x^2}-1 \\
& =&-\left(1+s^2\right) \\
X_1 & =&-\left(1+s^2\right) \frac{\partial}{\partial s}+s t \frac{\partial}{\partial t}\\
X_2 (s) & =&\left(\left(-\frac{y}{x}\right) \frac{\partial}{\partial x}+\left(x^2-\frac{y^2}{x^2}\right) \frac{\partial}{\partial y}\right) \frac{y}{x} \\
& =&\frac{y^2}{x^3}+x-\frac{y^2}{x^3}=t\\
X_2 (t)&=&\frac{-y}{x}=-s, \\
X_2&=&t\frac{ \partial}{\partial s}-s \frac{\partial}{\partial t}. \\
X_3 (s)&=&\left(\left(1+x^2\right) \frac{\partial}{\partial x}+\left(2 x y+\frac{y}{x}\right) \frac{\partial}{\partial y}\right) \frac{y}{x} \\
& =&-\left(1+x^2\right) \frac{y}{x^2}+2 y+\frac{y}{x^2} \\
& =&\frac{-y}{x^2}-y+2 y+\frac{y}{x^2}=y=st \\
 X_3 (t) &=&\left(1+x^2\right)=1+t^2,  \\
 X_3&=&s t \frac{\partial}{\partial s}+\left(1+t^2\right) \frac{\partial}{\partial t}. 
\end{eqnarray*}

Writing the realization of $\mathrm{SO}_3$ in last question
\begin{eqnarray*}
 X_1&=&y \frac{\partial}{\partial x}-x \frac{\partial}{\partial y}, \\
X_2&=&-x y \frac{\partial}{\partial x}+\left(-1-y^2\right) \frac{\partial}{\partial y}, \\
 X_3&=&(1+x^2) \frac{\partial}{\partial x}+x y \frac{\partial}{\partial y}.
\end{eqnarray*}
Comparing with this problem
\begin{eqnarray*}
X_1&=&-\left(1+s^2\right)\frac{\partial}{\partial s}+s t \frac{\partial}{\partial t}, \\
 X_2&=&t\frac{ \partial}{\partial s}-s \frac{\partial}{\partial t}, \\
X_3&=&s t+\left(1+t^2\right) \frac{\partial}{\partial t}.
\end{eqnarray*}

As the given algebra corresponds to the group of rotations, without integrating, the two first integrals could also be obtained alternatively.  Rotation is certainly a symmetry. Being a member of projective group it is going to throw lines to lines (solutions in st-plane). The other two symmetries are hard to visualize but can be checked by $X^{n} \omega= \eta^{n}$.\ \ $\triangleleft$ \\

If a chain of derived group terminates at the identity group then such a group is solvable. Easiest example of solvable groups are abelian groups. So $G_2$ is solvable. To show that let infinitesimal generators $X_1$ and $X_2$ be basis vectors for a two dimensional lie algebra. Suppose $\left[X_1,X_2\right]=k X_1+l X_2$ $=Z$. If $k$ and $l$ are zero then there is nothing to show as $L$ is abelian. If for $c_1$ and $c_2$ arbitrary constants; $c_1 X_1+c_2 X_2 \in L$  then
\begin{eqnarray*}
&& {\left[Z, c_1 X_1+c_2 X_2\right]=c_1\left[Z, X_1\right]+c_2\left[Z, X_2\right] }, \\
&&=c_1 l\left[X_2, X_1\right]+c_2 k\left[X_1, X_2\right] \\
&&=\left(c_2 k-c_1 l\right) Z
\end{eqnarray*}
hence $Z$ is an ideal of $L$ with dimension one. \ \ $\triangleleft$ \\
 
To study a connection between differential equations and solvability observe that a solvable $G_{n}$ is admitted by every linear differential equation of $n^{\text {th }}$ order. It is not hard to see that generators occurring as symmetries of linear differential equations are
$$X_N=U_N(x) \frac{\partial}{\partial y}.$$
Here $U_N$ are $n$ linearly independent solutions of the corresponding homogeneous linear differential equations. The symmetries form an abelain group as structure constants are all zero; $\left[X_N, X_M\right]=0$.\\
Note: Linear differential equations can have more than $n$ symmetries. Take $y^n=0$ for instance  has translation along $x$ axis and other symmetries as well.\ \ $\triangleleft$ \\

Lie's theory works for higher order differential equations even better. To see, solve $y^{\prime \prime \prime}=\frac{3}{2} \frac{y^{\prime \prime}}{y^{\prime}}$ using first strategy \cite{hs}. 
\begin{eqnarray*}
&& \text{Taking} \ X_1=\frac{\partial}{\partial x} ; X_2=x \frac{\partial}{\partial x}. \\
&& {\left[X_1, X_2\right]=\left[\frac{\partial}{\partial x}, x \frac{\partial}{\partial x}\right]=X_1}. \\
&&t=y, \ \ s=x. \\
&& X_1=\frac{\partial}{\partial s}, \\
&& X_2=x {\frac{\partial}{\partial x}}(y) \frac{\partial}{\partial x}+x \frac{\partial}{\partial x}(x) \frac{\partial}{\partial s}=s \frac{\partial}{\partial s}, \\
&& A=\frac{\partial}{\partial t}+ s^{\prime} \frac{\partial}{\partial s}+s^{\prime \prime} \frac{\partial}{\partial s^{\prime}}+(\cdots) \frac{\partial}{\partial s^{\prime \prime}}, \\
&& X_2^{(3)}=s\frac{ \partial}{\partial s}+s^{\prime} \frac{\partial}{\partial s^{\prime}}+s^{\prime \prime} \frac{\partial}{\partial s^{\prime \prime}}+s^{\prime \prime \prime} \frac{\partial}{\partial s^{\prime \prime \prime}}, \\
&& y^{\prime \prime \prime}=\frac{3 s^{\prime \prime^2}-s^{\prime} s^{\prime \prime \prime}}{s^{\prime^5}}, y^{\prime \prime}=\frac{-s^{\prime \prime}}{s^{\prime^3}}, \ \ y^{\prime}=\frac{1}{s^{\prime}}, \\
&& \frac{3 s^{\prime \prime^2}-s^{\prime} s^{\prime \prime\prime}}{s^{\prime^ 5}}=\frac{3}{2}\left(\frac{-s^{\prime \prime}}{s^{\prime^3}}\right)^2 s^{\prime}, \\
\end{eqnarray*}
\begin{equation} \label{s3prime}
 s^{\prime \prime \prime}=\frac{3}{2}\frac{s^{\prime\prime^2}}{s^{\prime}}.
\end{equation}
\begin{eqnarray*}
&& Y_{2}=s^{\prime}\frac{\partial}{\partial s^{\prime}}+s^{\prime \prime} \frac{\partial}{\partial s^{\prime\prime}} \\
&& Y_{2} \ u(t,s^{\prime})=1; \ \ Y_{2} \ v(t,s^\prime)=0\\
&& s^{\prime}\frac{\partial}{\partial s^{\prime}}+s^{\prime \prime} \frac{\partial}{\partial s^{\prime\prime}} \ (u)=1 \\
&& s^{\prime}\frac{\partial}{\partial s^{\prime}}+s^{\prime \prime} \frac{\partial}{\partial s^{\prime\prime}} \ (v)=0\\
&& u=\ln s^{\prime}, \ \  v=t \\
&& e^u=s^{\prime} \\
&& s^{\prime \prime}=e^u u^{\prime} \\
&& s^{\prime \prime \prime}=e^u u^{\prime \prime}+e^{u} v^{\prime 2} 
\end{eqnarray*}

(\ref{s3prime}) is transformed to
\begin{eqnarray*}
&&e^u u^{\prime}+e^{u} u^{\prime 2}=\frac{3}{2} \frac{\left(e^{u}{u^{\prime^2}}\right)^2}{ e^u} \\
&&u^{\prime \prime}+u^{\prime^2}=\frac{3}{2} u^{\prime^2} \\
&&u^{\prime \prime}=\frac{1}{2} u^{\prime^2}\\
\operatorname{let} &&u^{\prime}=w \\
&&u^{\prime \prime}=w^{\prime}\\
&&w^{\prime}  =\frac{1}{2} w^2 \\
&&\frac{d w}{w^2}  =\frac{1}{2} d v \\
&&-\frac{1}{w}  =\frac{1}{2} v+c_1 \\
&&-2  =v w+2 c_1 w \\
&&w  =\frac{-2}{2 c_1+v}\\
&&\frac{d u}{d v}  =\frac{-2}{c_1^{\prime}+v} \\
&&u  =-2 \ ln \mid c_1^{\prime}+v \mid+c_2 \\
&&ln \ s^{\prime} =ln\left|c_1^{\prime}+t\right|^{-2} c_2^{\prime}\\
&&s^{\prime} =\frac{c_2^{\prime}}{\left(c_1^{\prime}+t\right)^2} \\
&&\frac{d s}{d t}  =\frac{c_2^{\prime}}{\left(c_1^{\prime}+t\right)^2}\\
&& s=\frac{-c_2^{\prime}}{\left(c_1^{\prime}+t\right)}+c_3^{\prime} \\
&& x=\frac{-c_2^{\prime}}{\left(c_1^{\prime}+y\right)}+c_3^{\prime} \\
&& \left(c_1^{\prime}+y\right) x=-c_2^{\prime}+c_3 c_1^{\prime}+c_3^{\prime} y \\
&& -c_3^{\prime} y+y x=-c_1^{\prime}+c_3^{\prime} c_1^{\prime}-c_1^{\prime} x \\
&& y=\phi_1-\frac{\phi_2}{x-\phi_3}
\end{eqnarray*}

So the given differential equation is the equation of rectangular hyperbolas that includes all the isotherms. Notice that two straight quadrature could have solved the given differential equation  but goal was to implement first strategy.\ \ $\triangleleft$ \\

As a last example use  $X_1=\frac{\partial}{\partial x},
X_2=\frac{\partial}{\partial y}, X_3 = x\frac{\partial}{\partial x}+y \frac{\partial}{\partial y},
X_4=y \frac{\partial}{\partial x}-x \frac{\partial}{\partial y}$ to solve $y^{\prime \prime \prime}\left(1+y^{\prime^2}\right)=\left(3 y^{\prime^2}+a\right) y^{\prime \prime}$. If the group is not transitive then a linear relation between $A$ and the generators
$X_N$ exits, and the coefficients of the $X_N$ in the relation may lead to first integrals. One way is to take the matrix of operators, using Gaussian elimination and method of characteristics, find at least one differential invariant, this differential invariant will act as an orbit. But in the space of $\left(x, y, y^{\prime}, y^{\prime \prime}\right)$ the five operators $X_1, X_2, X_3, X_4$ and $A$ are always linearly dependent so
$$X_1=\phi_1 X_2+\phi_2 X_3+\phi_3 X_4+ \nu A$$
$$\frac{\partial}{\partial x}=\phi \frac{\partial}{\partial y}+\phi_2\left(x\frac{\partial}{\partial x}+y\frac{\partial}{\partial y}-y^{\prime\prime}\frac{\partial}{\partial y^{\prime\prime}}-2y^{\prime\prime\prime}\frac{\partial}{\partial y^{\prime\prime\prime}}\right)+ \\$$
$$\phi_3\left(y\frac{\partial}{\partial x}-x\frac{\partial}{\partial y}-(1+y^{\prime^2})\frac{\partial}{\partial y^{\prime}}-3y^{\prime}y^{\prime\prime}\frac{\partial}{\partial y^{\prime\prime}}-3y^{\prime\prime^2}-4y^{\prime}y^{\prime\prime\prime}\frac{\partial}{\partial y^{\prime\prime\prime}}\right)+\\$$
$$\nu \left(\frac{\partial}{\partial x}+y^{\prime}\frac{\partial}{\partial y}+y^{\prime\prime}\frac{\partial}{\partial y^{\prime}}+\omega\frac{\partial}{\partial y^{\prime\prime}}\right)$$
\begin{eqnarray*}
&& det=\left|\begin{array}{ccc}
x & y & 1 \\
0 & -\left(1+y^{\prime 2}\right) & y^{\prime \prime} \\
-y^{\prime \prime} & -3 y^{\prime} y^{\prime \prime} & \frac{\left(3 y^{\prime}+a\right) y^{\prime \prime}}{\left(1+y^{\prime 2}\right)}
\end{array}\right| \\
&& \ \ =x\left(-\left(3 y^{\prime}+a\right) y^{\prime \prime}+3 y^{\prime} y^{\prime\prime}\right)-y^{\prime \prime}\left(y y^{\prime \prime}+\left(1+y^2\right))\right. \\
&& a_{11}=\left|\begin{array}{lll}
y & -x & y^{\prime} \\
0 & -\left(1+y^{\prime^2}\right) & y^{\prime \prime} \\
y^{\prime \prime} & -3 y^{\prime} y^{\prime\prime} & \frac{\left(3 y^{\prime}+a\right) y^{\prime \prime^2}}{1+y^{\prime^2}}
\end{array}\right| \\
&& a_{11}=y\left(-\left(3 y^{\prime}+a\right) y^{\prime \prime^2}+3 y^{\prime} y^{\prime \prime^2}\right)-y^{\prime \prime}\left(-x y^{\prime \prime}+\left(1+y^{\prime^2}\right) y^{\prime}\right) \\
&& \phi_1=\frac{y\left(-\left(3 y^{\prime}+a\right) y^{\prime\prime^2}+3 y^{\prime} y^{\prime \prime^2}\right)-y^{\prime \prime}\left(-x y^{\prime \prime}+\left(1+y^{\prime^2}\right) y^{\prime}\right)}{x\left(-\left(3 y^{\prime}+a\right) y^{\prime \prime}+3 y^{\prime} y^{\prime \prime}\right)-y^{\prime \prime}\left(y y^{\prime \prime}+\left(1+y^{\prime^2}\right)\right.} \\
&& a_{12}=-\left|\begin{array}{ccc}
1 & -x & y^{\prime} \\
0 & -\left(1+y^{\prime^2}\right) & y^{\prime \prime} \\
0 & -3 y^{\prime} y^{\prime \prime} & \frac{\left(3 y^{\prime}+a\right) y^{\prime \prime^2}}{1+y^{\prime^2}}
\end{array}\right| \\
&& \phi_2=\frac{-\left(-\left(3  y^{\prime}+a\right) y^{\prime \prime^2}+3 y^{\prime} y^{\prime \prime^2}\right)}{x\left(-\left(3 y^{\prime}+a\right) y^{\prime \prime}+3 y^{\prime} y^{\prime \prime}\right)-y^{\prime \prime}\left(y y^{\prime \prime}+\left(1+y^{\prime^2}\right)\right.} \\
&&a_{13}=\left|\begin{array}{ccc}1 & y & y^{\prime} \\ 0 & 0 & y^{\prime \prime} \\ 0 & -y^{\prime \prime} & \frac{\left(3 y^{\prime}+a\right) y^{\prime\prime^2}}{\left(1+y^{\prime 2}\right)}\end{array}\right|\\
&&\phi_3=\frac{y^{\prime \prime^2}}{x\left(-\left( 3y^{\prime}+a\right) y^{\prime\prime^2}+3 y^{\prime} y^{\prime \prime^2}-y^{\prime \prime}\left(y y^{\prime \prime}+\left(1+y^{\prime^2}\right)\right)\right)}  
\end{eqnarray*}\ \ $\triangleleft$ \\
This concludes our modest cookbook on how to solve odes using symmetry analysis. For pde see \cite{blusym} 
\section{Linearizing FHND}
Now that isolated fixed points have been found for CFHND. Stability analysis for solutions near fix points can be carried out easily. To perform the stability analysis one needs to derive the Jacobian for FHND first.  Recall \cite{rou} that stability for an equilibrium solution for a first order ordinary differential equation $\frac{d u}{d t}=f(u)$  can be studied doing the following procedure. Firstly perturbing $u$ about a fix point $u_{*}$. \\

$u(t)=u_{*}+\delta(u)(t),$ \\ 

$\frac{d}{d t}\left(u_{*}+\delta u\right)=f\left(u_{*}+\delta u\right),$ \\ 

$\frac{d}{d t}\left(u_{*}+\delta u\right)=\frac{d \delta u}{d t},$ \\ 

$f\left(u_{*}+\delta u\right) \approx f\left(u_{*}\right)+f^{\prime}\left(u_{*}\right) \delta u,$ \\ 

$\frac{d \delta u}{d t}=f^{\prime}\left(u_{*}\right),$ \\ 

$  \delta u=a e^{\sigma t},$ \\  

$\sigma =  f^{\prime} (u_*),$ \\  

$ \delta u(t) = a \ e^{{f^{\prime}(u_*)t}}, $  \\

Now if $f^{\prime}\left(u_{*}\right)<0 \ \ u_{*}$ is linearly stable \\  \\
If $f^{\prime}\left(u_{*}\right)>0 \ \ u_{*}$ is linearly unstable \\  \\
Next given reaction diffusion equation $\partial_{t} u(t, x)=D \partial_{x x}^{2} u+f(u),$ stability of an isolated fix point can be  seen repeating a similar idea of perturbing substituting and solving.\\

$u(t, x)=u_{*}+\operatorname{\delta u}(t, x),$ \\  

$\partial_{t}\left(u_{*}+\delta u(t, x)\right)= \partial_{x x}^{2}\left(u_{*}+\delta u\right)+f (u_{*}+\delta u), $ \\  

$\partial_{t} \operatorname{\delta u}(t, x)=D \partial_{x x}^{2} \delta u+f^{\prime}\left(u_{*}\right) \delta u, $ \\  

\begin{eqnarray*} \delta u(t, x) &&=a e^{\sigma t} e^{i k x}, \\ 
\sigma &&=f\left(u_{k}\right)-D k^{2}.  \end{eqnarray*} 
Now for a system of ordinary differential equations $$\frac{d u}{d t}=f_{1}(u, v);
\  \frac{d v}{d t}=f_{2}(u, v),$$ the analysis is not different. \\

$u(t)= u_{*}+ \delta u(t),$ \\

$v(t)=v_{*}+\delta v(t);$ \\  

$\frac{d}{d t} [u_{*}+\delta u(t)]= \frac{d \delta u}{d t}= f_{1} (u_{*}+\delta u , v_{*}+\delta v ),$ \\ 

$\frac{d}{d t}[v_{*}+\delta v(t)]=\frac{d \delta v}{d t}=f_{2}(u_{*}+\delta u, v_{*}+\delta v),$ \\ 

$\frac{d \delta u}{d t}=a_{11} \delta u+a_{12} \delta v,$ \\ 

$\frac{d \delta v}{d t}=a_{21} \delta u+a_{22}  \delta v,$ \\ 

$a_{i j}=\frac{\partial f_{i}}{\partial u_{j}}\mid_{\left(u^{*},v^{*}\right)},$ \\ 

$\delta u(t)=\delta u_{0} e^{\sigma t},$ \\

$\delta v(t) = \delta v_{0} e^{\sigma t};$ \\

$ \sigma \delta u_{0}=a_{11} \delta u_{0}+a_{12} \delta v_{0},$ \\

$\sigma \delta v_{0}=a_{21} \delta u_{0}+a_{22}  \delta v_{0},$ \\ 

$A \delta \vec{u}_{0}=\sigma \delta \vec{u}_{0}.$ \\ 

\noindent So finally let us gather these ideas to linearize FHND about an isolated fix point. \\

$\partial_{t} u(t, x)=D \partial_{x x}^{u} u+f_{1}(u, v),$ \\

 $\partial_{t} v(t, x)= f_{2}(u,v);$ \\ 

$u(t, x)=u_{*}+\delta u(t, x),$ \\

$v(t, x)=v_{*}+\delta v(t, x);$ \\ 

$\partial_{t} \delta u(t, x)=\ \ D \partial_{xx}^{2} \delta u + a_{11} \delta u + a_{12} \delta v,$ \\

$\partial_{t} \delta v(t, x)=  a_{21} \delta u+a_{22} \delta v;$ \\ 

$ \delta u(t, x)= \delta u_{0} e^{\sigma t} e^{i k x},$ \\

$\delta v(t, x)=\delta v_{0} e^{\sigma t} e^{i k x};$ \\

 $J=\left(\begin{array}{cc}a_{11}-D k^{2} & a_{12} \\ a_{21} & a_{22}\end{array}\right).$

\section{Non-classical symmetry analysis}
Qualitative analysis for differential equations mostly become meaningless if exact solution to the system is known. Lie symmetry analysis is one of the most powerful tool available to solve differential equations exactly and systematically.  Especially for nonlinear differential equation where no general methods are available.\\
(\ref{fd})  has two obvious Lie symmetries. The spatio-temporal translation invariance. These two symmetries do not  give any differential invariant so they are not enough to run quadrature. G. Bluman \cite{Bluman1967}  extended the idea of Lie symmetries. He introduced this idea of non-classical symmetries. Most of the times it is difficult to tackle a system of differential equations.
The system of two equations (\ref{fd}) were merged together and as a result a third order differential equation was obtained:\\

\begin{equation}\label{3pd}
\mathrm{D} u_{t x x}-u_{t t}+\varepsilon  \beta  \mathrm{D} u_{x x}-\varepsilon  \beta  u_{t}-\varepsilon  u + u_{t} \left(-u^{2}+1\right)+ 
 \varepsilon  \beta  \left(u -\frac{u^{3}}{3}\right) =0\;.
\end{equation}

\section{Infinite Exact Solutions of Fitzhugh-Nagumo Model with Diffusion}
When non-classical symmetry analysis (Appendix) for conditional symmetries \cite{FMNAG} were carried out on (\ref{3pd}), following non-classical  symmetries were obtained.
\begin{eqnarray}
 Y_{1}:& \xi_{x}=1\\
 Y_{2}:& \xi_{t}=1, \xi_{x}= \sqrt{-\mathrm{D} \beta \epsilon}\\
Y_{3}:& \xi_{t}=1, \xi_{x}= - \sqrt{-\mathrm{D} \beta \epsilon}\\
Y_{4}:& \xi_{t}=1, \xi_{x}= c\\
Y_{\infty}:& \xi_{t}=1, \eta_{u}= A(t,x) u+ B(t,x) 
\end{eqnarray}

Using $Y_{\infty}$ following ansatze (\ref{ansa}) together with the constraints (\ref{constraint1}-\ref{constraint last}) resulted. 
\begin{eqnarray}   
&u(x,t)=\frac{{\mathrm e}^{A t} \textit{F}(x) A-B}{A}\label{ansa}\\
&-\epsilon \beta \textit{F}(x)^{3} A^{3}-3 \textit{F}(x)^{3} A^{4}=0\label{constraint1}\\
&3 \epsilon \beta \textit{F}(x)^{2} B \,A^{2}+6 \textit{F}(x)^{2} B \,A^{3}=0\\
&3 \epsilon \beta \mathrm{D} (\frac{d^{2}}{dx^{2}}\textit{F}(x)) A^{3}+3 \mathrm{D} (\frac{d^{2}}{dx^{2}}\textit{F}(x)) A^{4}-3 \epsilon \beta \textit{F}(x) A^{4} \nonumber \\
&-3 A^{5} \textit{F}(x)+3 \epsilon \beta \textit{F}(x) A^{3}-3 \epsilon \beta \textit{F}(x) B^{2} A  +3 \textit{F}(x) A^{4} \nonumber\\& -3 \epsilon \textit{F}(x) A^{3}-3 \textit{F}(x) B^{2} A^{2}=0\\
&-3 \epsilon \beta B \,A^{2}+\epsilon \beta B^{3}+3 \epsilon B \,A^{2}=0\label{constraint last}
\end{eqnarray}

These differential and algebraic constraints produced following results.

\begin{eqnarray}
\textit{F}(x)&=&\mathit{c_{1}} {\mathrm e}^{\frac{\sqrt{A^{3} \beta \epsilon+A^{4}-A^{2} \beta \epsilon+B^{2} \beta \epsilon-A^{3}+A^{2} \epsilon+A \,B^{2}}\, x}{A \sqrt{\mathrm{D}}\, \sqrt{\epsilon \beta+A}}}\nonumber\\
&+&\mathit{c_{2}} {\mathrm e}^{-\frac{\sqrt{A^{3} \beta \epsilon+A^{4}-A^{2} \beta \epsilon+B^{2} \beta \epsilon-A^{3}+A^{2} \epsilon+A \,B^{2}}\, x}{A \sqrt{\mathrm{D}}\, \sqrt{\epsilon \beta+A}}}
\end{eqnarray}
\begin{equation}
A=-\frac{\epsilon \beta}{3}, B=0  
\end{equation}
\begin{equation} \label{ncsolu}
u(t,x)={\mathrm e}^{-\frac{\epsilon \beta t}{3}} \left(\mathit{c_{1}} {\mathrm e}^{-\frac{\sqrt{-2 \epsilon^{4} \beta^{4}-6 \epsilon^{3} \beta^{3}+9 \epsilon^{3} \beta^{2}}\, x \sqrt{6}}{6 \epsilon \beta \sqrt{\mathrm{D}}\, \sqrt{\epsilon \beta}}}+\mathit{c_{2}} {\mathrm e}^{\frac{\sqrt{-2 \epsilon^{4} \beta^{4}-6 \epsilon^{3} \beta^{3}+9 \epsilon^{3} \beta^{2}}\, x \sqrt{6}}{6 \epsilon \beta \sqrt{\mathrm{D}}\, \sqrt{\epsilon \beta}}}\right)
\end{equation}

\begin{equation}
k=\frac{\sqrt{-2 \epsilon^{4} \beta^{4}-6 \epsilon^{3} \beta^{3}+9 \epsilon^{3} \beta^{2}}\,  \sqrt{6}}{6 \epsilon \beta \sqrt{\mathrm{D}}\, \sqrt{\epsilon \beta}}
\end{equation}

\begin{equation}
u(t,x)={\mathrm e}^{-\frac{\epsilon \beta t}{3}} (\mathit{c_{1}} {\mathrm e}^{-kx}+\mathit{c_{2}} {\mathrm e}^{kx})
\end{equation}
\begin{eqnarray}\label{ncsolv}
v(t,x)=-\frac{1}{6 \beta ({\mathrm e}^{\frac{\epsilon \beta t}{3}})^{3} ({\mathrm e}^{kx})^{3}}[2 ({\mathrm e}^{kx})^{6}(c_{2})^{3} \beta+6 ({\mathrm e}^{kx})^{4} c_{1}(c_{2})^{2} \beta\nonumber\\
-9 ({\mathrm e}^{kx})^{4} ({\mathrm e}^{\frac{\epsilon \beta t}{3}})^{2} c_{2}+6 ({\mathrm e}^{kx})^{2}(c_{1})^{2} c_{2} \beta-9 ({\mathrm e}^{kx})^{2} ({\mathrm e}^{\frac{\epsilon \beta t}{3}})^{2} c_{1}+2(c_{1})^{3}\beta]
\end{eqnarray}

\section{A particular solution of FitzHugh-Nagumo}

\begin{figure}[h!] 
\centering
\includegraphics[scale=0.5]{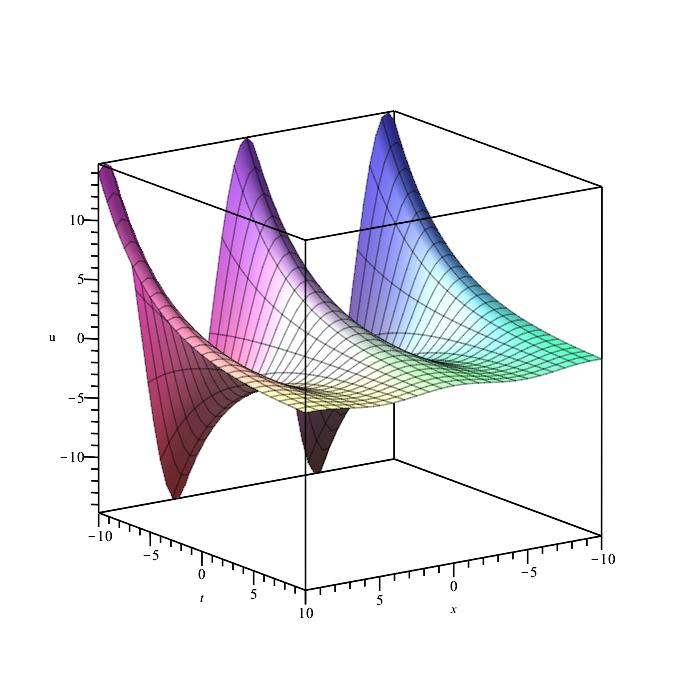}
\caption{An exact solution (\ref{ncsolu}) for FHND (\ref{fd}); A solution surface $u(t,x)$;  $\mathit{c_{1}}=1,\mathit{c_{2}}=1,\epsilon=0.3,\beta=2,D=1.03,c=0$}
\label{ncl}
\end{figure}

(\ref{ncsolu}) and (\ref{ncsolv}) to the best of our knowledge is the first ever exact solution to (\ref{fd}). Although the existence of traveling wave (\ref{ncsolu}) was indicated (Fig.\ref{ncl}) by numerous numerical and asymptotical works done earlier \cite{rinzel},\cite{gns}. The greatest gift symmetry analysis has presented to the age-long FHND is the algorithm from $Y_{\infty}$ to generate infinite exact solutions for FHND. \\

\begin{equation}
u(t,x)={\mathrm e}^{-\frac{\epsilon \beta t}{3}} (\mathit{c_{1}} {\mathrm e}^{-kx}+\mathit{c_{2}} {\mathrm e}^{kx})
\end{equation}

\begin{equation}
k=\frac{\sqrt{-2 \epsilon^{4} \beta^{4}-6 \epsilon^{3} \beta^{3}+9 \epsilon^{3} \beta^{2}}\,  \sqrt{6}}{6 \epsilon \beta \sqrt{\mathrm{D}}\, \sqrt{\epsilon \beta}}
\end{equation}
\begin{figure}[h!] 
\centering
\includegraphics[scale=0.5]{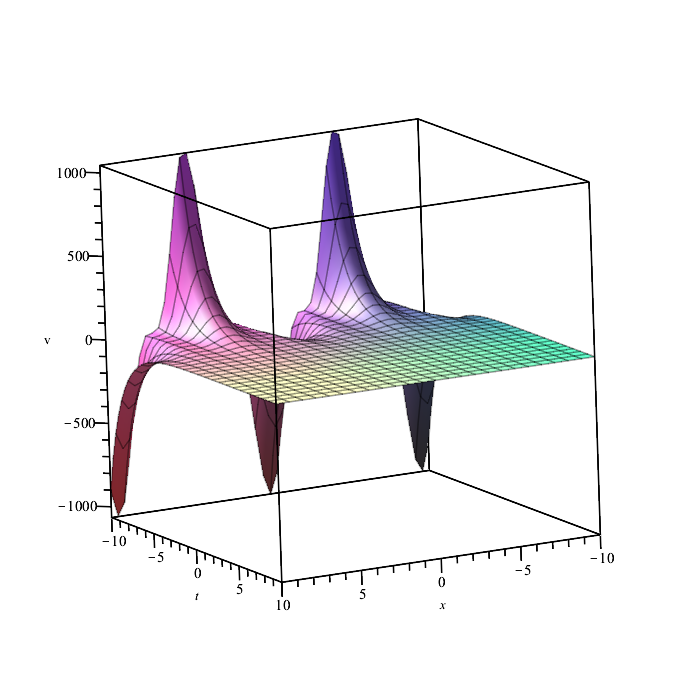}
\caption{An exact solution (\ref{ncsolv}) for FHND (\ref{fd}); A solution surface $v(t,x)$;  $\mathit{c_{1}}=1,\mathit{c_{2}}=1,\epsilon=0.3,\beta=2,D=1.03,c=0$}
\label{vncl}
\end{figure}

\section{Appendix}
Following is partial calculation for non-classical symmetries of (\ref{fd}). Consider  
\begin{eqnarray*}
u_t &=& D u_{xx} - v + g(u)  \\
v_t &=& \epsilon \left[ -\beta v + c + u \right]  
\end{eqnarray*}
Writing differential consequence of the fast equation
$$ u_{tt} = D u_{xxt} - v_t + g_{t}(u) $$
 
Substituting the value of $v_t$
$$  u_{tt} = D u_{xxt} - \epsilon \left[ -\beta v + c + u \right] + g_t(u)    $$
Making $v$ the subject in the first fast FHND
$$ v = D u_{xx} - u_t + g(u)  $$
Substituting this value in the differential consequence
$$u_{tt} - D u_{xxt} + \epsilon \left[ -\beta (D u_{xx} - u_t + g(u)) \right]  + c + u ] - g_t(u) = 0 $$
FHND has now been transformed into the following third order pde
$$ D u_{xxt} - u_{tt} + \epsilon \beta D u_{xx} + g_t(u) - \epsilon \beta u_t + \epsilon \beta g(u) - \epsilon c - \epsilon u = 0 $$
Consider next the infinitesimal generator for non-classical symmetry and its extension 
$$ \Gamma = \xi^o \frac{\partial}{\partial t} + \xi^1 \frac{\partial}{\partial x} + \eta \frac{\partial}{\partial u}$$ 
And its third extension
$$\Gamma^{(3)} = \Gamma + \rho_x \frac{\partial}{\partial u_x} + \rho_t \frac{\partial}{\partial u_t} + \sigma_{tt} \frac{\partial}{\partial u_{tt}} 
    + \sigma_{tx} \frac{\partial}{\partial u_{tx}} + \sigma_{xx} \frac{\partial}{\partial u_{xx}} + \tau_{xxt} \frac{\partial}{\partial u_{xxt}} + \ldots$$

 $$D_t \tau_{xxt} - \sigma_{tt} + \epsilon \beta D \sigma_{xx} + \eta g_{tu} - \epsilon \beta \rho_t + \epsilon \beta \eta g_u - \eta \epsilon = 0 $$
    $$\eta, \rho_t, \sigma_{xx}, \sigma_{tt}, \tau_{xxt} $$
    $$\rho_t = D_t (\eta) - u_t D_t (\xi^o) - u_x D_t (\xi^1) $$
    $$D_t = \frac{\partial}{\partial t} + u_t \frac{\partial}{\partial u} + u_{tt} \frac{\partial}{\partial u_t} + u_{tx} \frac{\partial}{\partial u_x}$$
    $$\rho_t = \eta_t + u_t \eta_u - u_t \xi_t^o - u_t^2 \xi_t^o - u_x \xi^1_t - u_{t}u_x \xi^1_u $$
    $$\sigma_{tt} = D_t (\rho_t) - u_{tt} D_t (\xi^o) - u_{tx} D_t (\xi^1) $$
   \begin{eqnarray*}
        &\sigma_{tt} = \eta_{tt} + u_t \eta_{ut} - u_t \xi_t - u_t \xi_t - 2 u_t \xi_{tt} - 2 u_t \xi_t - u_x \xi_{tt} + \eta u_{tt} \\ 
    &- u_{tt} \xi_u - u_{tx} \xi_x - u_x \xi_{xt} + u_t \eta_{tt} - u_x \xi_t - u_t \xi_x - u_{tx} \xi_u + u_{tx} \xi - u_{tx} \xi - u_{tt} \xi \\
    &- u_x^2 \xi_u - u_x^2 \xi_t - u_{tt} \xi_x - u_x \xi_t - 2 u_t \xi_{tx} - u_{tx} \eta_x - u_x^2 \eta - u_t \xi - u_{tx} \eta_{tx} + \eta_t - u_t \xi_u \\
    &- u_x^2 \eta_u - 2 u_t \xi_u - u_t \xi_x - u_{tx} u_x \xi - 2 u_{tx} u_x \xi_u - 2 u_{tx} u_x \xi_x - u_{tx} \eta_{tx} - u_{tx} \xi_t - u_t u_x \xi_u \\
    &- u_t u_x \xi_t - 2 u_{tx} \eta_x - u_t \xi_{tx} - 2 u_t \xi_t - u_x^2 \xi_u - 2 u_{tx} \xi_x - u_x \eta_u - u_x \xi_t - u_{tx} \xi_u - u_{tx} \xi \\
    &- u_{tx} \eta_u - u_x^2 \xi_t - u_x \xi_t - u_x \xi_x - u_t \xi - u_{tx} \eta_x - u_x \xi_t - u_x \xi_t - u_t \xi_{tt} - u_{tx} \xi_t - u_x \eta_u \\
    &- u_t \eta_x - u_{tx} \xi_t - u_x \xi_t - u_x^2 \xi - u_{tx} \xi_t - u_{tx} \xi_u 
    \end{eqnarray*}
    \\ \\
    $$\sigma_{xx} = D_x(\rho_{x}) - u_{tx} D_x(\xi^o) - u_{xx} D_x(\xi^1) $$
    \begin{eqnarray*}
        \rho_x &=& D_x(\eta) - u_t D_x(\xi^o) - u_x D_x(\xi^1) \\
    &=& \eta_x + u_x \eta_u - u_t \xi^o_x - u_t u_x \xi_u^o - u_x \xi_x^1 - u_x^2 \xi_u^1 
    \end{eqnarray*}
\begin{eqnarray*}  
    \sigma_{xx} &=& \eta_{xx} + u_x \eta_{ux} - u_t \xi_{xx}^o - u_t u_x \xi_{ux}^o \\
    &\quad& - u_x \xi_{xx}^1 - u_x^2 \xi_{ux}^1 + u_x [\eta_u + u_x \eta_{ux} - u_t \xi_x^o - u_x \xi_u^1] \\
    &\quad& - 2 u_x \xi_{ux} - u_{tx} [\xi_x^o + u_x \xi_u^o] - u_{xx} [\xi_x^1 + u_x \xi_u^1]
\end{eqnarray*}
\begin{eqnarray*}
    \tau_{xxt} &=& D_t(\sigma_{xx}) - u_{tx} D_t(\xi^0) - u_{xxt} D_t(\xi^1) \nonumber \\
    &=& D_x(\sigma_{tx}) - u_{tx} D_x(\xi^0) - u_{xxt} D_x(\xi^1) \nonumber \\
    &=& \eta_{xxt} + u_t \eta_{uxu} + u_x \eta_{uxu} + u_t u_x \eta_{uxu} + u_x u_x \eta_{uxu} \nonumber \\
    &\quad &- u_t \xi_{xxt} - u_t u_x \xi_{uxxt} - u_t^2 \xi_{uxt} - u_t u_{xx} \xi_{uxu} - u_t u_x \xi_{uxt} \nonumber \\
    &\quad &- u_t u_x \xi_{uxu} - u_t u_{xx} \xi_{uxu} - u_x u_t \xi_{uxu} - u_x u_t \xi_{uxu} \nonumber \\
    &\quad &- u_{xxt} [\xi^0 + u_t \xi_u^0] - u_{xxt} [\xi^1 + u_x \xi_u^1] - u_{tx}^2 \xi_{xx}^1 - u_t^2 \xi_{uu} \nonumber \\
    &\quad &- u_t u_x \xi_{uu} - u_t \xi_u^1 - u_t \xi_u^1 - 2 u_t u_{xx} \xi_{uu} + u_x \eta_{xu} \nonumber \\
    &\quad& + u_t u_x \eta_{ut} + u_{xx} \eta_{uu} + u_{tx} \eta_{ut} + u_x \eta_{uu} + u_x^2 \eta_{uu} \nonumber \\
    &\quad& - u_{tx} u_x \xi_{ux} + u_{tx} \eta_u + \eta_{tx} \xi_u - u_{tx} \xi_{tx} + u_t \eta_{uu} \nonumber \\
    &\quad& + \eta_{tx} \eta_u + \eta_{tx}^2 \xi_{tx} + u_x \eta_{xx} + u_t \eta_{tx} + u_t \eta_{ut} - u_{tt} \eta_u \nonumber \\
    &\quad& - u_{tx} \eta_t + u_x^2 \eta_{uu} + u_t u_{tx} \eta_u + u_{tx} \eta_u + \eta_{tx} \eta_u - u_{tx} \xi_{tx} \nonumber \\
    &\quad& + u_t \eta_{ut} - u_t u_x \xi_{uu} + u_{xx} \eta_{uu} + u_x^2 \eta_{ut} - u_x^3 \xi_{ut} + u_x \eta_{xx} \nonumber \\
    &\quad& + \eta_{tx} \xi_{tt} + u_x^2 \eta_{uu} - 2 u_t u_{xx} \xi_{tx} + u_x \eta_{xx} + u_t u_{xx} \eta_{ut} \nonumber \\
    &\quad& + u_{xx} \eta_{ut} + \eta_{tx} \xi_u - 2 u_t u_{tx} \xi_{uu} - u_{tx} u_x \xi_{tx} + u_x^2 \eta_{tx} \nonumber \\
    &\quad& + u_{tx} \eta_u - u_{tx} \xi_{tx} - u_{tx} u_x \eta_{ut} + u_t u_{tx} \xi_{ut} + u_x^2 \eta_{uu} + u_t \eta_{tx} \nonumber \\
    &\quad& - u_{tt} \eta_u - 2 u_t u_x \xi_{uu} - u_x u_{xx} \eta_{ut} + \eta_{tx} \eta_{tx} + u_t^2 \eta_{ut} \nonumber \\
    &\quad& + u_x \eta_{tx} + u_t u_{xx} \eta_{ut} + u_x^2 \eta_{tx} - u_x u_{xx} \eta_{ut} + u_{tx} \eta_u + u_x \eta_{xx} \nonumber \\
    &\quad& - u_{tx} u_x \eta_{ut} + u_{tx} \eta_u + u_{tx} \eta_u - u_t u_{xx} \xi_{ut} + u_t \eta_{xx} + u_{tx} \eta_{ut} \nonumber \\
    &\quad& + u_{tx} \eta_u - u_{tx} u_x \xi_{ut} - u_x u_{xx} \eta_{ut} + \eta_{tx} \xi_u - u_t u_x \xi_{ut} - u_x^2 \eta_{tx} \nonumber \\
    &\quad& + u_t u_{tx} \eta_u + u_t^2 \eta_{tx} + u_x \eta_{xx} - u_t u_{xx} \eta_{ut} + u_x \eta_{xx} - u_t u_x \xi_{ut} \nonumber \\
    &\quad& - u_{tx} u_x \eta_{ut} + \eta_{tx} \xi_u + u_t^2 \eta_{ut} + u_{tx} \eta_u + u_x^2 \eta_{tx} + u_x \eta_{xx} \nonumber \\
    &\quad& - u_{tx} \xi_{tx} - u_x u_{xx} \eta_{ut} - u_{xx} \xi_{xt}^1 - u_t u_x \xi_{xt}^1 - u_{tx} \xi_{xt}^1 - u_x u_x \xi_{xt}^1 \nonumber \\
    &\quad& - u_t u_x u_{xt} \xi_{ut}^0 - u_x u_x \xi_t^1 - u_{tx} u_x \xi_u^1 - u_{txx} \xi_u^1 - u_{tx} \xi_t^1 - u_t u_{txx} \xi_u^1 \nonumber \\
    &\quad& + D u_{xxt} + D u_x \nu_{xt} + D u_x \nu_{xt} + D u_x \nu_{xt} + D u_x \nu_{xt} + D u_{tx} \nu_{xt} \nonumber \\
    &\quad& - D u_t \xi_{ut}^0 - D u_t \xi_{ut}^0 - D u_t \xi_{ut}^0 - D u_x u_{tx} \xi_{ut}^0 - D u_x \xi_{ut}^0 - D u_{tx} u_t \xi_u^1 \nonumber \\
    &\quad& - D u_{tx} \xi_{ut}^1 - D u_x \xi_{ut}^1 - D u_{tx} u_x \xi_t^1 - D u_t \xi_u^1 - D u_x \xi_{xt}^1 - D u_t u_x \xi_{xt}^1 \nonumber \\
    &\quad& - D u_t^2 \xi_u^1 - D u_t u_{xt}^1 \xi_{ut} - 2 D u_t u_{tx} \xi_u^1 + D u_x \nu_{xt} + D u_{xt} \nu_{xt} \nonumber \\
    &\quad& + D u_{tx} u_{xt} + D u_t \xi_{xt} + D u_{xt} u_x + 2 D u_{tx} u_{xt} \nu_{xt} - D u_x \xi_{ut}^1 \nonumber \\
    &\quad& - D u_{xt} \xi_{ut}^1 - D u_{tx} \xi_u^1 - D u_t u_{tx} \xi_{ut}^0 - D u_{tx} u_x - D u_{tx} \xi_{ut} - D u_t u_{tx} \xi_{ut}^1 \nonumber \\
    &\quad& - D u_t u_x \xi_{ut} - D u_{tx} u_t \xi_{ut} - 2 D u_{tx} u_x \xi_t - D u_{tx} \xi_t - D u_t u_{tx} \xi_t \nonumber \\
    &\quad& - D u_{tx} u_{tx} \xi_u^1 - D u_t u_{tx} - D u_t u_x \xi_{xt} - D u_{tx} \xi_u - D u_t u_{tx} \xi_{ut} \nonumber \\
    &\quad& - D u_{tx} u_t \xi_{ut}^1 - D u_{tx} u_x \xi_{ut} - D u_{tx} \xi_t^0 - D u_{tx} \xi_x - \eta_t - u_t \nu_t + u_t \xi_{tt}^1 \nonumber \\
    &\quad& + u_t \xi_{tt}^1 - D u_{ttx} \xi_x^0 - D u_x u_{tx} \xi_{ut}^0 - D u_t u_{tx} \xi_{xt}^1 - D u_{tx}^2 \xi_t^1 \nonumber \\
    &\quad& - D u_{ttx} \xi_t^1 - D u_{ttx} u_x \xi_n^1 + D u_{txx} u_t + D u_{tx} u_{tt} u_t \nonumber \\
    &\quad& + D u_{tx} \nu_t - D u_{tx} u_t \xi_t^1 - D u_{txx} \xi_t^1 - D u_{ttx} u_x \xi_{xt}^1 - D u_{tx} u_x \xi_{xt}^1 \nonumber \\
    &\quad& - D u_{tx} \xi_{xt}^1 - 2 D u_{tx} u_x \xi_{xt}^1 - D u_{txx} \xi_t^1 - 2 D u_{tx} u_x \xi_t^1 - 2 D u_{ttx} u_x \xi_{xt}^1 \nonumber \\
    &\quad& - D u_{txx} \xi_{xt}^1 - D u_{tx} \xi_t^1 - D u_{ttx} \xi_t^1 - D u_{txx} \xi_{xt}^1 \nonumber \\
    &\quad& - D u_{txx} \xi_t^1 - D u_{tx} u_x \xi_{xt}^1 - D u_{ttx} u_x \xi_t^1 - D u_{tx} \xi_t^1 - D u_{txx} u_x \xi_t^1 \nonumber \\
    &\quad& - D u_{tx} \xi_{xt}^1 - D u_{txx} u_t - D u_{tx} \xi_t^1 - D u_{txx} \xi_t^1 - D u_{tx} \xi_{xt}^1 \nonumber \\
    &\quad& - D u_{tx} u_x \xi_{xt}^1 - D u_{txx} u_x \xi_t^1 - D u_{tx} \xi_{xt}^1 - D u_{txx} \xi_{xt}^1 \nonumber \\
    &\quad& - D u_{tx} \xi_t^1 - D u_{txx} \xi_t^1 - D u_{tx} \xi_{xt}^1 - D u_{txx} u_x \xi_{xt}^1 \nonumber \\
    &\quad& - D u_{tx} u_x \xi_{xt}^1 - D u_{tx} \xi_t^1 - D u_{txx} \xi_t^1 - D u_{tx} \xi_{xt}^1 - D u_{txx} \xi_t^1 \nonumber \\
    &\quad& - D u_{tx} \xi_{xt}^1 - D u_{txx} u_x \xi_{xt}^1 - D u_{tx} u_x \xi_{xt}^1 - D u_{tx} \xi_t^1 - D u_{txx} \xi_t^1 \nonumber \\
    &\quad& - D u_{tx} \xi_{xt}^1 - D u_{txx} \xi_t^1 - \eta_t - u_t \nu_t + u_t \xi_t^0 + u_t \xi_{tt}^1 \nonumber \\
    &\quad& + u_t u_x \xi_{ut}^1 - u_t \eta_u - u_{tt} \eta_u + u_x^2 \xi_t + u_t^2 \xi_u + u_t u_x \xi_u + u_x^2 \xi_{uu} \nonumber \\
    &\quad& - \eta_u u_{tt} + u_{tt} \xi_t + u_{tt} u_x \xi_u + u_t \xi_t + u_t u_{xx} \xi_{tt} + \epsilon \beta D \eta_{ux} + \epsilon \beta \eta_u \nonumber \\
    &\quad& - \epsilon \beta D u_t \xi_u u_x - \epsilon \beta D u_t \xi_{ux} - \epsilon \beta D u_t^2 \xi_t u_{xx} - \epsilon \beta D u_x \xi_t \nonumber \\
    &\quad& - \epsilon \beta D u_x^2 \xi_{tt} - \epsilon \beta D u_{xx} u_x \xi_t - \epsilon \beta D u_t^2 \xi_{tt} u_x - \epsilon \beta D u_t \xi_t \nonumber \\
    &\quad& - \epsilon \beta D u_x \xi_{tt} - \epsilon \beta D u_x u_{xx} \xi_u - \epsilon \beta D u_x u_{tt} \xi_t - \epsilon \beta D u_t \xi_x \nonumber \\
    &\quad& - \epsilon \beta D u_x \xi_{tt} - \epsilon \beta D u_x u_t \xi_u - \epsilon \beta D u_x u_{xx} \xi_t - \epsilon \beta D u_x \xi_u \nonumber \\
    &\quad& + \epsilon \beta D u_x \xi_u + \epsilon \beta D u_x u_t \xi_t + \epsilon \beta D u_x \xi_{tt} + \epsilon \beta D u_x u_x \xi_t \nonumber \\
    &\quad& + \epsilon \beta D u_t \xi_t + \epsilon \beta u_x \xi_u + \epsilon \beta D u_x \xi_u - \epsilon \eta = 0
\end{eqnarray*}

\begin{eqnarray*}
 \text{Const}&: & \quad D \eta_{xxt} -  \eta_{tt} + \epsilon \beta D \eta_{xuu} + \xi_{uu} - u_{xx} \xi \\
                  & \quad &- \epsilon \beta \xi_{ut} + \epsilon \beta \eta_{uu} - \epsilon \eta = 0 \\[10pt]
    u_{t}&: & \quad D \eta_{uxu} - D \xi_{xuu} - \eta_{uuu} + \xi_{tt} - \eta_{ut} + \xi \\
           & \quad& - 2 \eta_{uu} + \xi_{ut} - u_{x}^2 \eta + u_{t}^2 \xi + \xi_{t} - \epsilon \beta D \xi_{xxu} \\
           & \quad& - \epsilon \beta \eta = 0 \\[10pt]
    u_{x}&: & \quad D \eta_{xuu} - D \xi_{xuu} + D \eta_{xxu} + \xi_{t} + \xi_{t} u_{x}^2 \\
           & \quad& + \epsilon \beta \eta_{uu} - \epsilon \beta \xi_{xuu} + \epsilon \beta \xi_{ut} = 0 \\[10pt]
    u_{t} u_{x}&: & \quad D \eta_{xuu} - D \xi_{uxu} - \eta_{xuu} \\
                  & \quad& - \xi_{ut} + \epsilon \beta \eta_{xuu} - \epsilon \beta \xi_{u} = 0 \\[10pt]
    u_{x}^2&: & \quad D \eta_{uuu} - D \eta_{uxx} + D \eta_{xxu} - D \eta_{uxx} \\
              & \quad& - \xi_{ux} + \epsilon \beta \xi_{uxu} - \epsilon \beta \xi_{uxu} = 0 \\[10pt]
    u_{t} u_{x}^2&: & \quad - \xi_{ux} + \xi_{ux} u_{t}^2 + \epsilon \beta \eta_{uu} - D \xi_{xuu} \\
                    & \quad& - \eta_{uuu} + \epsilon \beta \xi_{uu} - \epsilon \beta D \xi_{u} = 0 \\[10pt]
    u_{tt}&: & \quad D \xi_{uxx} - \eta_{uu} + \xi_{t} + \xi_{t} = 0 \\[10pt]
    u_{t}^2 u_{x}&: & \quad \xi_{uu} - D \xi_{uxu} - D \xi_{xuuu} - D \xi_{uxuu} = 0 \\[10pt]
u_{tt} u_{x}&: & \quad - D \xi_{ux}^o - D \xi_{ux} + \xi_{uu} = 0 \\[10pt]
    u_{tx} u_t&: & \quad - D \xi_{ux}^o - D \xi_{uxu} - D \xi_{uxu} = 0 \\[10pt]
    u_{x}^2&: & \quad D \eta_{uut} - 2 D \xi_{xuu} - \epsilon \beta D \xi_{ux}^1 \\
             & \quad& + \epsilon \beta D \eta_{uuu} - \epsilon \beta D \xi_{xu} = 0 \\[10pt]
    u_{t} u_{x}^2&: & \quad - D \xi_{ux}^1 + D \eta_{uu} - D \xi_{uut} \\
                   & \quad& - D \xi_{xuu} - \epsilon \beta D \xi_{uu}^o = 0 \\[10pt]
    u_{tx} u_{xx}&: & \quad - 2 D \xi_{ux}^1 + 2 D \eta_{uuu} - 2 D \xi_{xuu} \\
                   & \quad& - D \xi_{ut} - \epsilon \beta D \xi_{u} - \epsilon \beta D \xi_{u} = 0 \\[10pt]
    u_{x}^2&: & \quad - D \xi_{uuu}^o = 0 \\[10pt]
    u_{t} u_{x}^2&: & \quad - D \xi_{uu}^o = 0 \\[10pt]
    u_{x}^3&: & \quad - D \xi_{uu}^1 - \epsilon \beta D \xi_{uux} = 0 \\[10pt]
    u_{t} u_{x}^3&: & \quad - D \xi_{uuu} = 0\\[10pt]
u_x^2 u_{tx}&: & \quad\quad\quad - 3 D \xi_{uu}^1 = 0 \\[10pt]
    u_t u_x&: & \quad\quad\quad - D \xi_{x}^o - D \xi_{x}^o = 0 \\[10pt]
    u_t u_{xx} u_{tx}&: & \quad\quad\quad - 2 D \xi_{uu}^1 - D \xi_{xu}^o - D \xi_{uu}^o = 0 \\[10pt]
    u_t^2 u&: & \quad\quad\quad - D \xi_u^o - D \xi_u^o = 0 \\[10pt]
    u_{tx} u_x&: & \quad\quad\quad - D \xi_u^o - D \xi_u^o = 0 \\[10pt]
    u_{xx}&: & \quad\quad\quad D \eta_{tt} - D \xi_{xt}^1 - D \xi_{xt}^1 + \epsilon \beta D \eta_u - \epsilon 3 D \xi_x^1 - \epsilon \beta D \xi_{u}^1 = 0 \\[10pt]
    u_t u_{xxt}&: & \quad\quad\quad - \epsilon \beta D \xi_u^o + D \eta_{uu} - D \xi_{ut}^o - D \xi_{xu}^1 - D \xi_{xu}^1 = 0 \\[10pt]
    u_t^2 u_{xx}&: & \quad\quad\quad - D \xi_{uu}^o = 0 \\[10pt]
    u_{xx}&: & \quad\quad\quad D \eta_{u} - D \xi_{x}^1 - D \xi_{x}^1 - D \xi_{t}^o = 0\\[10pt]
    u_{tx} u_{xx}&: & \quad\quad\quad - D \xi_u^o = 0 \\[10pt]
    u_{txx} u_t&: & \quad\quad\quad - D \xi_u^o - D \xi_u^o = 0 \\[10pt]
    u_{xx} u_x u_x&: & \quad\quad\quad - 2 D \xi_{ut}^1 - D \xi_{ut}^1 - 2 \epsilon \beta D \xi_u^1 = 0 \\[10pt]
    u_t u_x u_{xx} u_{xx}&: & \quad\quad\quad - 2 D \xi_{uu}^1 - D \xi_{uu}^1 - 0 = 0 \\[10pt]
    u_{xx} u_{xt}&: & \quad\quad\quad - D \xi_u^1 = 0 \\[10pt]
    u_{tx} u_{xx}&: & \quad\quad\quad - 2 D \xi_u^1 - D \xi_u^1 = 0 \\[10pt]
    u_{xx} u_{xt}&: & \quad\quad\quad - D \xi_t^1 = 0 \\[10pt]
    u_{xxt} u_t&: & \quad\quad\quad - D \xi_u^1 = 0 \\[10pt]
    u_t u_{tt}&: & \quad\quad\quad + 2 \xi_u^o + \xi_u^o = 0 \\[10pt]
    u_t^3&: & \quad\quad\quad \xi_{uu}^o = 0
\end{eqnarray*}
    
\begin{eqnarray*}
      \xi^o &=& E^o(t,x,u) \\
    &=& E^o(t,x) \\
    &=& E^o(t) \\
    \xi^1 &=& E^1(t,x,u) \\
    &=& E^1(t,x) \\
    &=& E^1(x) \\
    D \eta_{xxt} - \eta_{tt} + \epsilon \beta D \eta_{xx} + \eta g_{tu} + \epsilon \beta \eta_{t} 
    + \epsilon \beta \eta g_u - \epsilon \eta &=& 0 \\
    D \eta_{xxu} - \eta_{ut} + \xi_{tt}^o - \eta_{tu} - \epsilon \beta \eta_u 
    + \epsilon \beta \xi_t^o &=& 0 \\
    2 D \eta_{uxt} + 2 \epsilon \beta D \eta_{ux} - \epsilon \beta D \xi_{xx}^1 &=& 0 \\
    D \eta_{ux} - D \xi_{xx}^1 + D \eta_{xu} &=& 0 \\
    - \eta_{u} + 2 \xi_t^o &=& 0 \\
    \eta_{uu} &=& 0 
\end{eqnarray*}\cite{sad}

\includepdf[pages=1-2]{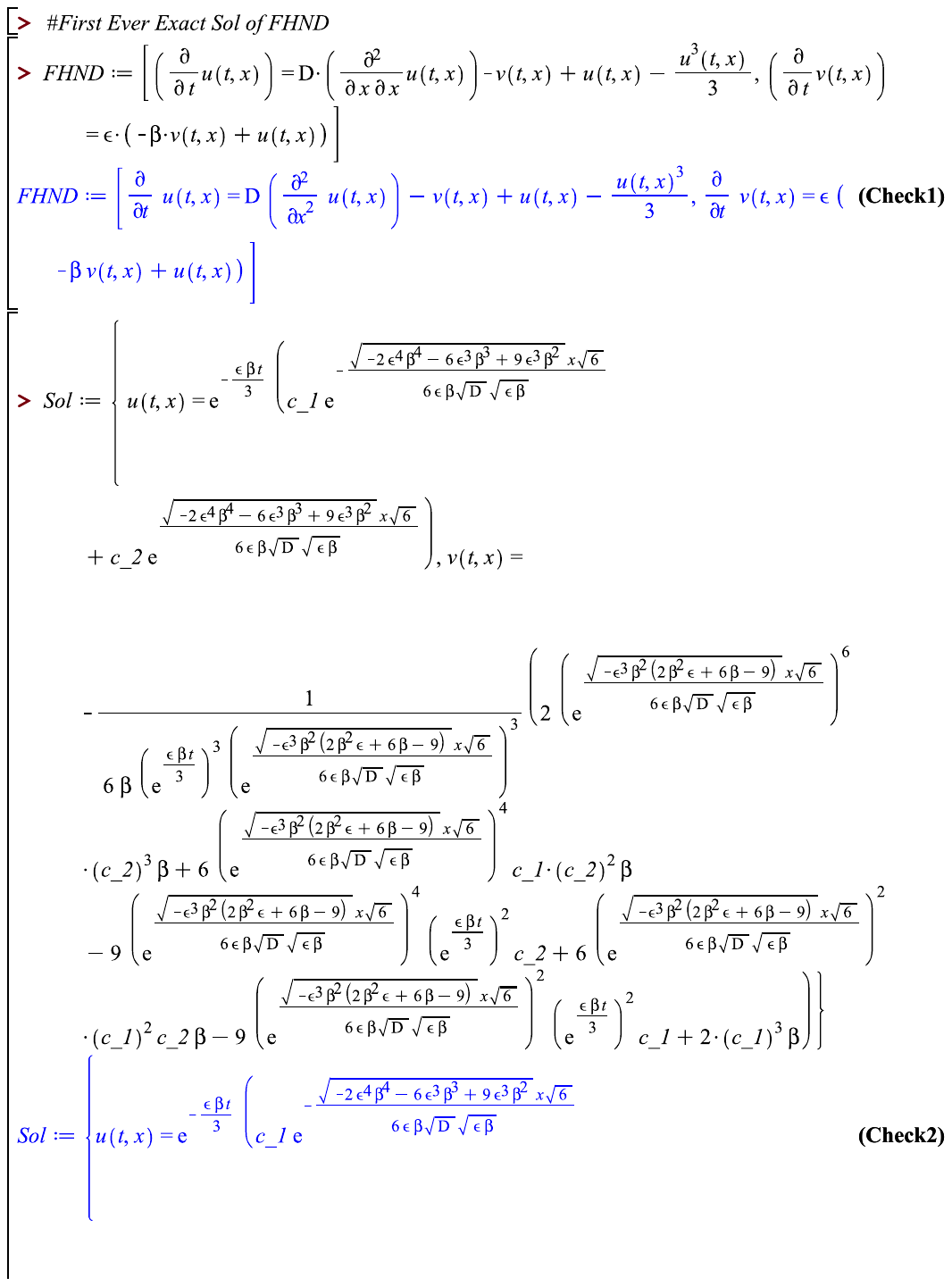}
\includepdf[pages=1-2]{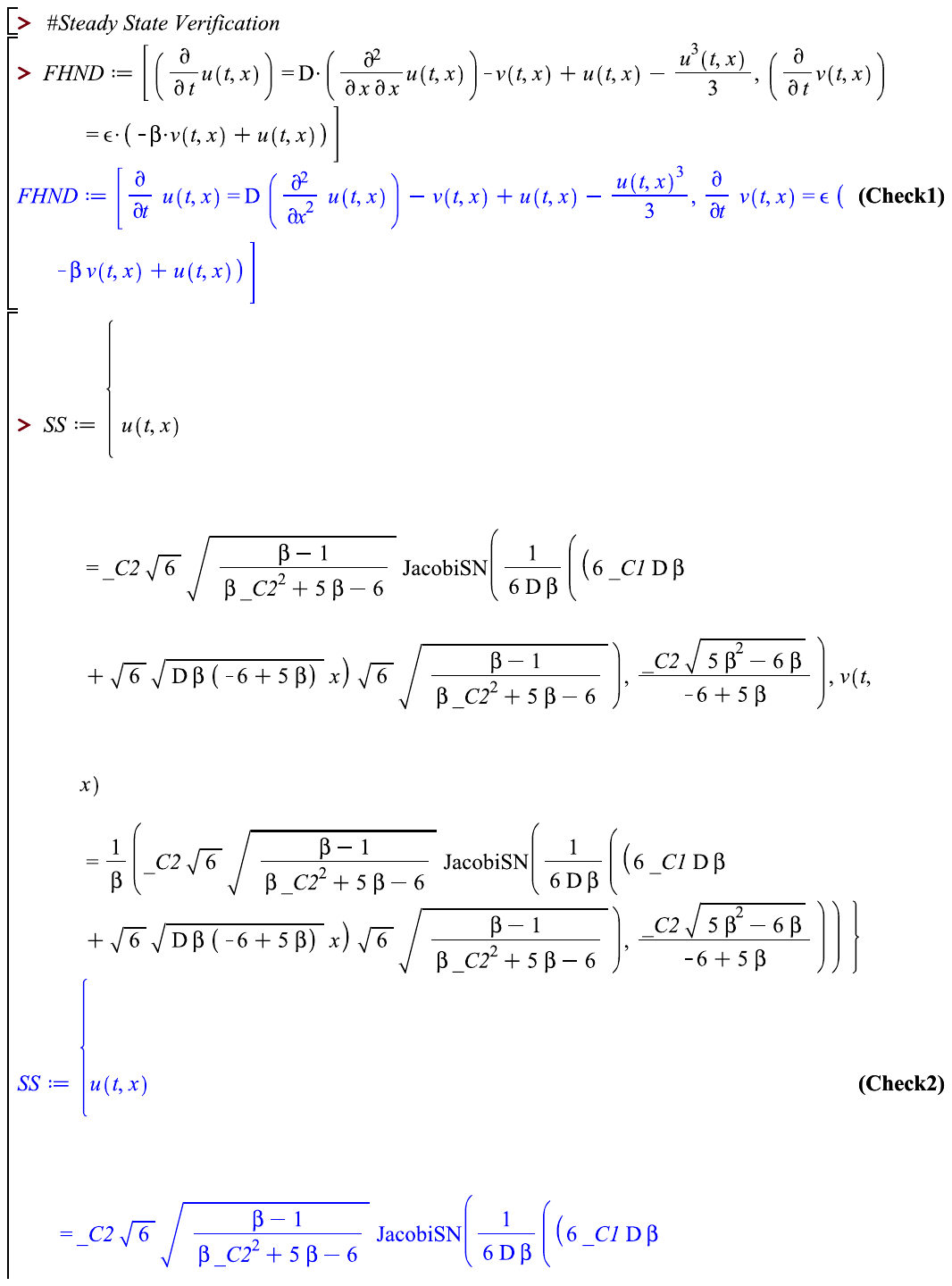}
\includepdf{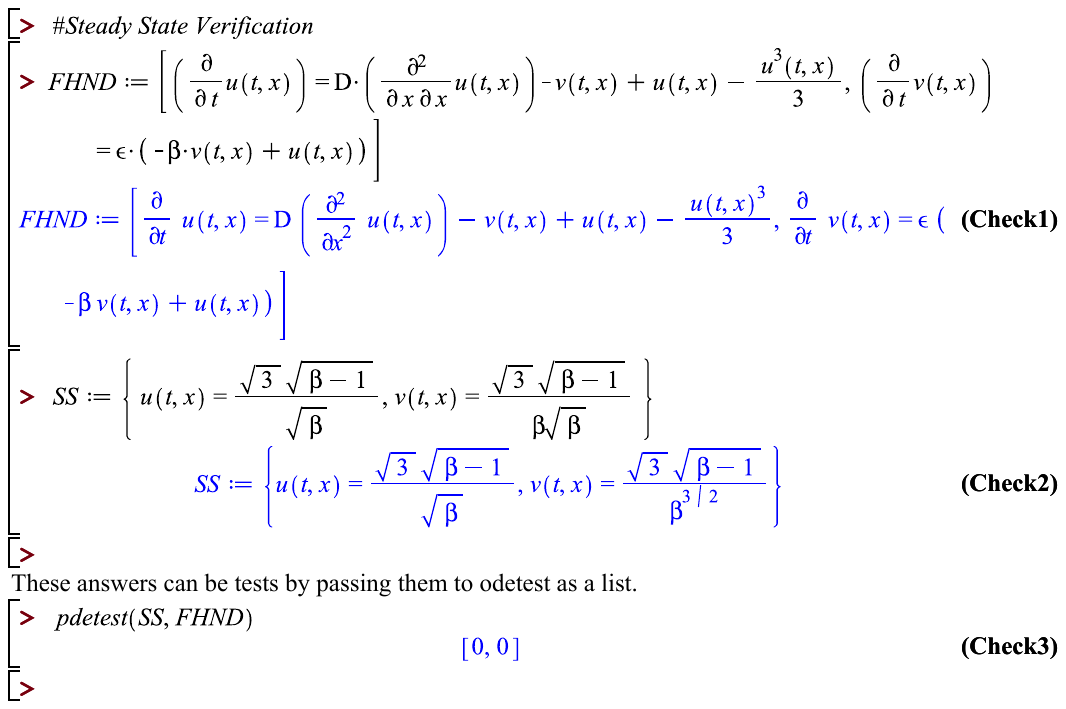}
\includepdf{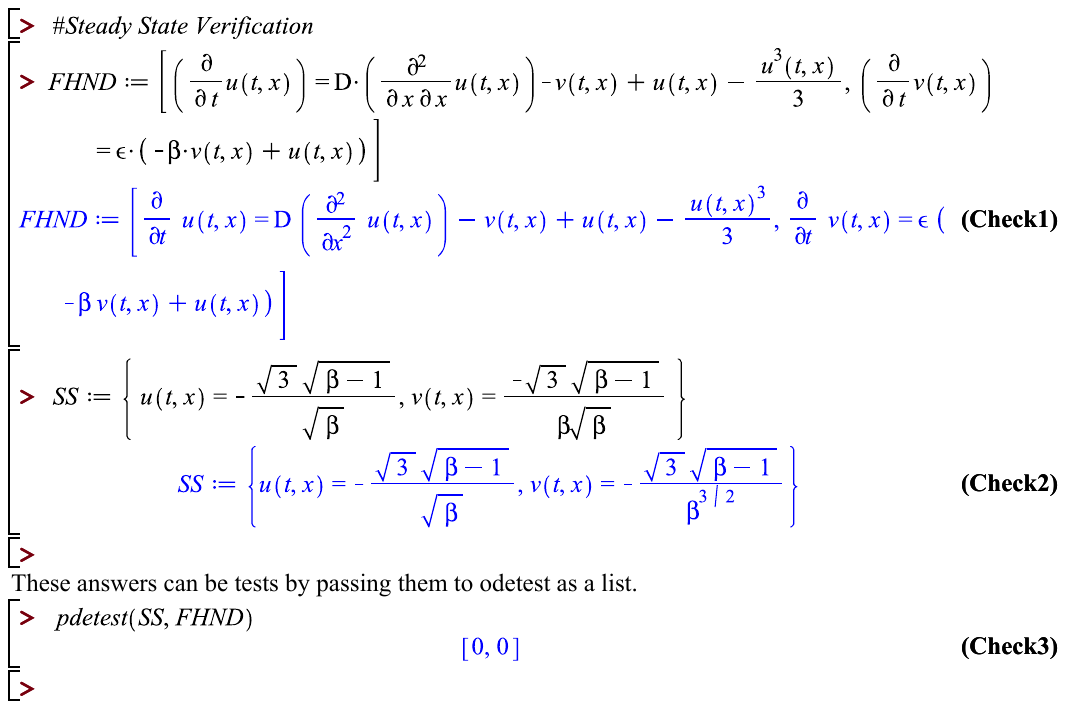}
\includepdf{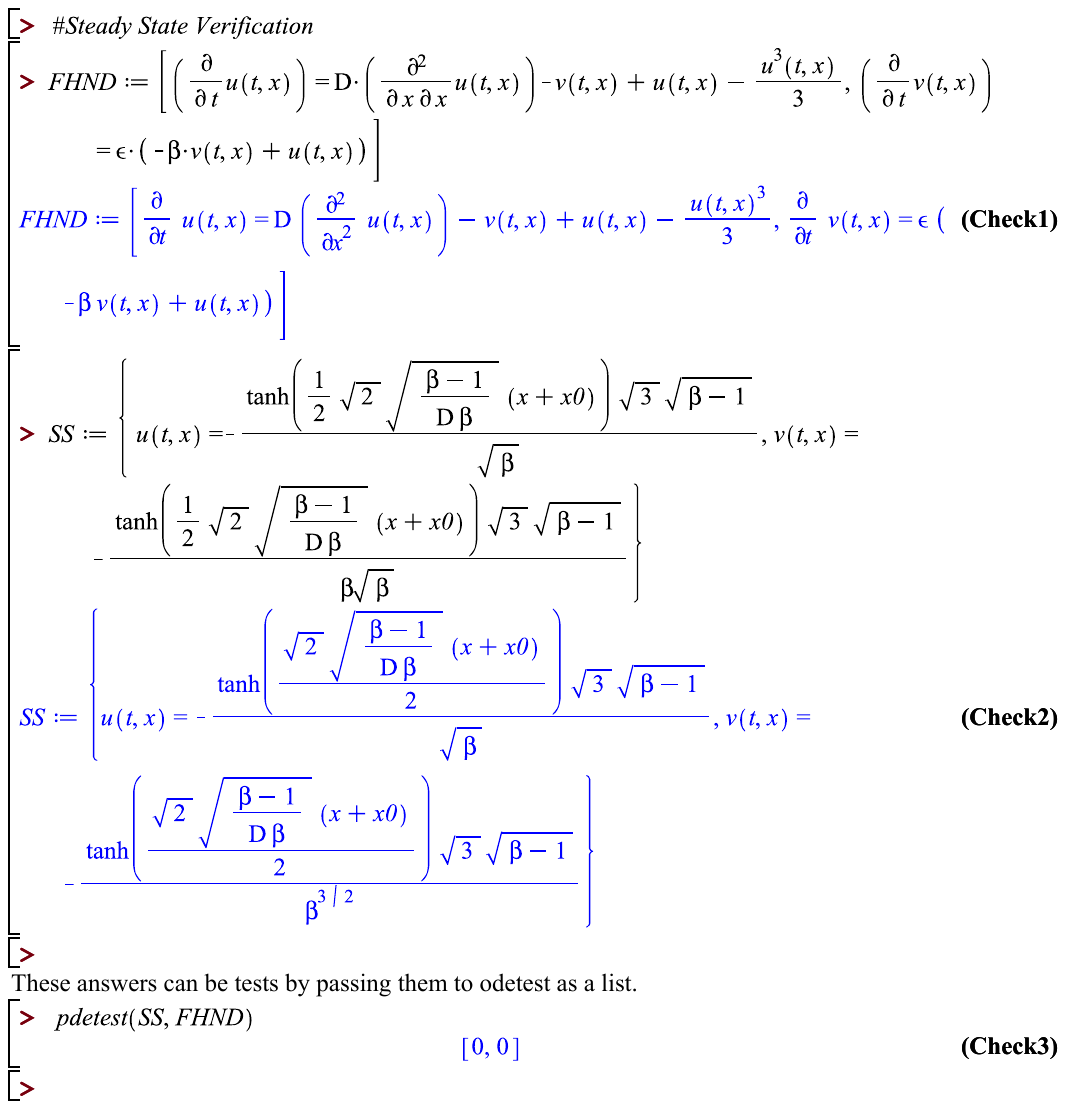}
\includepdf{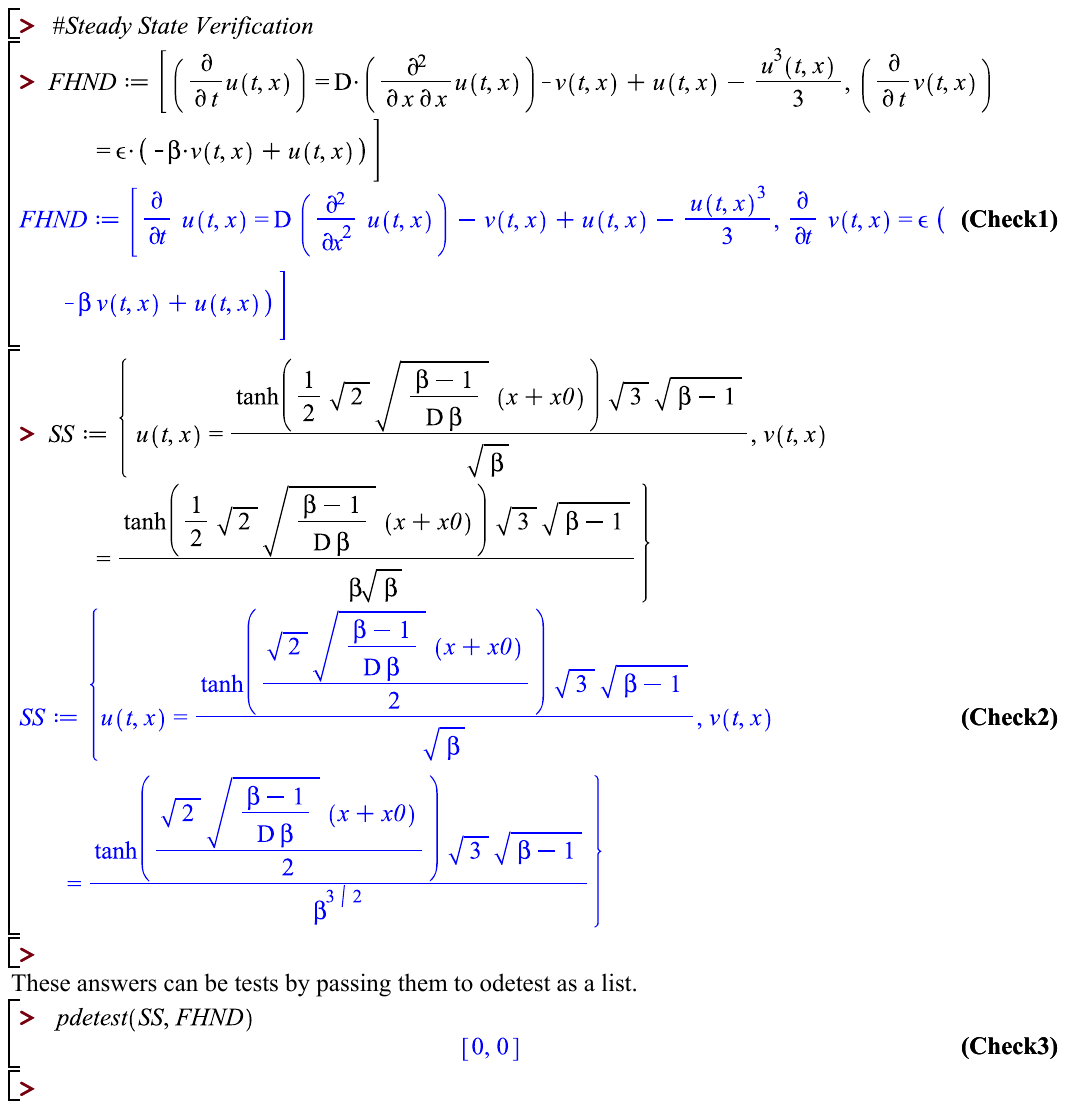}
\includepdf[pages=1-2]{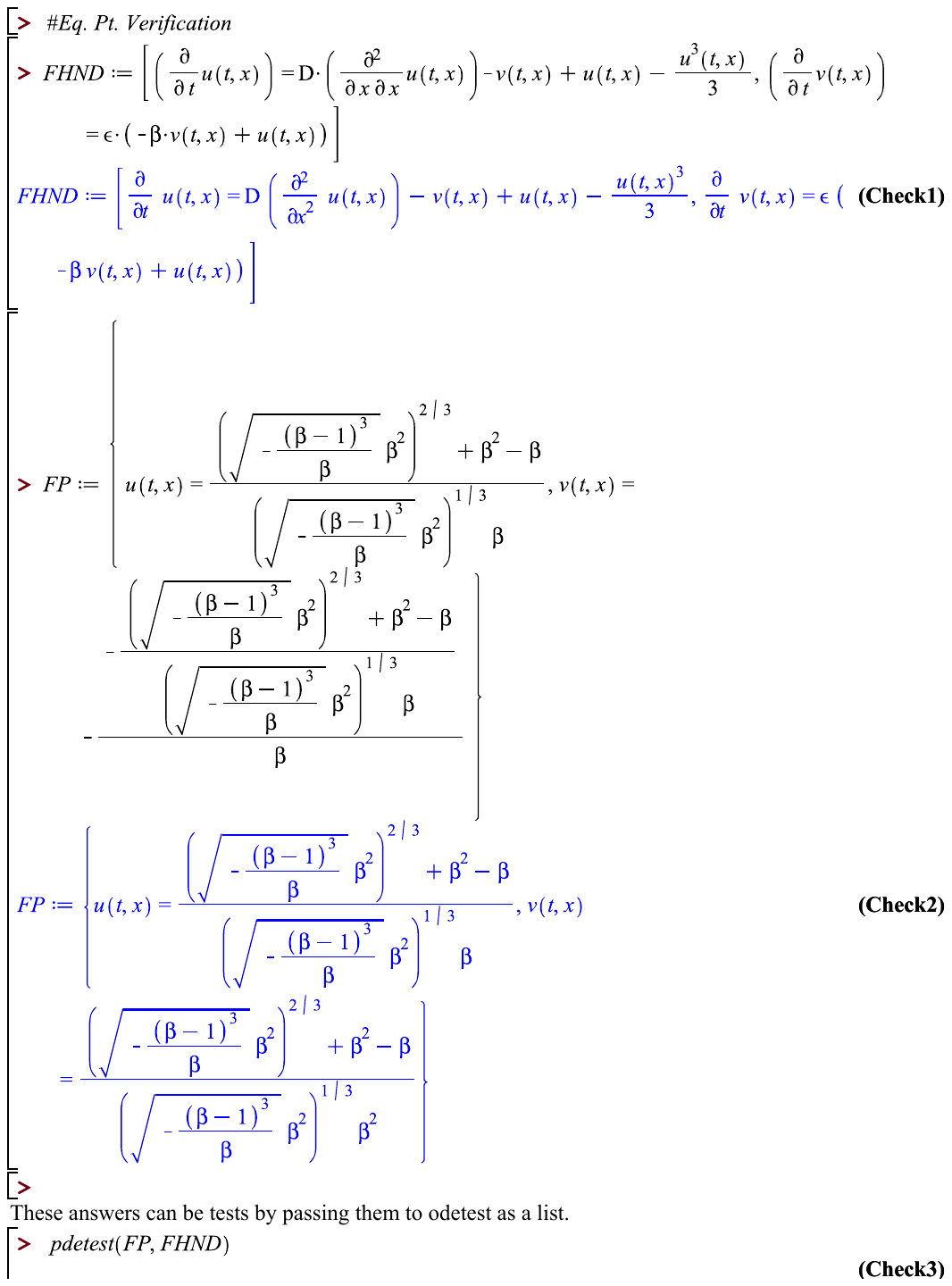}



\section*{Resource Availability Statement}
The following resources are available on request from the first author: Hand written calculations for non-classical symmetry for FitzHugh-Nagumo with diffusion; Matlab code for FitzHugh-Nagumo with/without diffusion; Hand calculations of related system using symmetry analysis \cite{zheng}; Conditional symmetry calculations for locally and non-locally related systems of FitzHugh-Nagumo with diffusion \cite{Cherniha2012}; Maple code for exp and Ricatti methods of traveling wave solutions; and Remaining steady state solutions of FitzHugh- Nagumo model.

\section*{Acknowledgements}
The first author expresses his heartfelt gratitude to Prof. Dr. Asghar Qadir, who proselytized him into symmetries.

%

 \section*{Conflict of interest}
The authors declare that they have no conflict of interest.


	%


\end{document}